\documentclass{Class-ICCM}

\usepackage{amsmath,amssymb,amsthm,mathrsfs,epsf,a4,color,graphicx,eucal}
  \usepackage{paralist}
\usepackage{graphics} %% add this and next lines if pictures in esp format
\usepackage{epsfig} %For pictures: screened artwork set up with an 85 or 100 line screen
\usepackage[colorlinks=true]{hyperref}
\hypersetup{urlcolor=blue, citecolor=red}
\usepackage{import}

\theoremstyle{plain}

\theoremstyle{definition}
\newtheorem{definition}{Definition}
\newtheorem{proposition}{Proposition}

\newtheorem{remark}{Remark}

\usepackage{amssymb}

\usepackage{psfrag}
\newcounter{myfigure}
\newenvironment{my-picture}[3]{\refstepcounter{myfigure}\label{#3}\setlength{\unitlength}{\textwidth}\begin{picture}(#1,#2)}{\end{picture}}
\usepackage{mathrsfs}   % loading \mathscr

\newcommand\calE{\mathscr E}
\newcommand\calR{\mathscr R}

\newcommand{\GD}{\mathchoice
                  {\Gamma_{\hspace*{-.15em}\mbox{\tiny\rm D}}}
                  {\Gamma_{\hspace*{-.15em}\mbox{\tiny\rm D}}}
                  {\Gamma_{\hspace*{-.1em}\mbox{\tiny\rm D}}}
                  {\Gamma_{\hspace*{-.05em}\mbox{\tiny\rm D}}}}
\newcommand{\GN}{\Gamma_{\hspace*{-.15em}\mbox{\tiny\rm N}}}
\newcommand\GC{\mathchoice{\Gamma_{\hspace*{-.15em}\mbox{\tiny\rm C}}}
                          {\Gamma_{\hspace*{-.15em}\mbox{\tiny\rm C}}}
                          {\Gamma_{\hspace*{-.05em}\mbox{\tiny\rm C}}}
                          {\Gamma_{\hspace*{-.05em}\mbox{\tiny\rm C}}}}

\newcommand\SC{\Sigma_{\mbox{\tiny\rm C}}}
\newcommand\uD{u_{\text{\tiny\rm D}}}

\newcommand\wD{w_{\text{\tiny\rm D}}}
\newcommand\DTwD{\DT w_{\text{\tiny\rm D}}}
\newcommand\R{\mathbb R}
\newcommand\N{\mathbb N}
\newcommand\pl{\partial}
\renewcommand\d{\mathrm d}
\newcommand\Diss{\mathrm{Diss}}
\newcommand\JUMP[2]{\mathchoice
                  {\big[\hspace*{-.3em}\big[#1\big]\hspace*{-.3em}\big]_{#2}}
                   {[\hspace*{-.15em}[#1]\hspace*{-.15em}]_{#2}}
                   {[\![#1]\!]_{#2}}
                   {[\![#1]\!]_{#2}}}
\newcommand\DT[1]{\mathchoice
                 {{\buildrel{\hspace*{.1em}\text{\LARGE.}}\over{#1}}}
                 {{\buildrel{\hspace*{.1em}\text{\Large.}}\over{#1}}}
                 {{\buildrel{\hspace*{.1em}\text{\large.}}\over{#1}}}
                 {{\buildrel{\hspace*{.1em}\text{\large.}}\over{#1}}}}

\newcommand{\DDD}[3]{\begin{array}[t]{c}#1\vspace*{-1em}\\_{#2}\vspace*{-.3em}\\_{#3}\end{array}}
\newcommand{\ddd}[3]{\DDD{\begin{array}[t]{c}\underbrace{#1}\vspace*{.6em}\end{array}}{\text{\footnotesize #2}}{\text{\footnotesize #3}}}

\newcommand{\weak}{\rightharpoonup}
\newcommand{\weaks}{\stackrel{*}{\rightharpoonup}}
\newcommand{\nablaS}{\nabla_{\scriptscriptstyle\textrm{\hspace*{-.2em}S}}}

\begin{document}
\begin{frontmatter}

%\markboth{T.~Roub\'\i\v cek, C.G.~Panagiotopoulos, V.~Manti\v c}{Local-solution
%approach to quasistatic rate-independent mixed-mode delamination}

\title{Local-solution approach to quasistatic
rate-independent mixed-mode delamination}

\author{Tom\'a\v s Roub\'\i\v cek}
\address{Mathematical Institute, Charles University,
Sokolovsk\'a 83, CZ-186~75~Praha~8,  Czech Republic\\ and \\
Institute of Thermomechanics, Czech Academy of Sciences,\\
Dolej\v skova 5, CZ-182~00~Praha~8, Czech Republic.\\
           roubicek@karlin.mff.cuni.cz       
}

\author{Christos G.~Panagiotopoulos}

\address{Institute of Applied and Computational Mathematics, FORTH,
Nikolaou Plastira 100,\\Vassilika Vouton,
GR-700 13 Heraklion, Crete, Greece, pchr@iacm.forth.gr}

\author{Vladislav Manti\v c}

\address{Department of Continuum Mechanics, School of Engineering,
University of Seville,\\[-.3em]
Camino de los Descubrimientos s/n, ES-41092 Seville, Spain,
mantic@us.es}

%\maketitle

%\begin{history}
%\received{(Day Month Year)}
%\revised{(Day Month Year)}
%\accepted{(Day Month Year)}
%\comby{(xxxxxxxxxx)}
%\end{history}

\begin{abstract}
The quasistatic rate-independent evolution of a delamination at small
strains in the so-called mixed mode, i.e.~distinguishing opening (Mode I)
from shearing (Mode II), devised in \cite{RoKrZe11DACM,RoMaPa13QMMD}, is
rigorously analyzed in the context of a concept of stress-driven local
solutions. The model has separately convex stored energy and is associative,
namely the 1-homogeneous potential of dissipative force driving the
delamination depends only on rates of internal parameters.
An efficient  fractional-step-type 
semi-implicit discretisation in time is shown to converge to
(specific,  stress-driven like) local solutions that may approximately 
obey the maximum-dissipation principle.
Making still a spatial discretisation,
this convergence as well as relevancy of such solution concept
are demonstrated on a nontrivial 2-dimensional example.
\end{abstract}

\end{frontmatter}
%\ccode{AMS Subject Classification: 
%35K85, %Unilateral problems and variational inequalities for parabolic PDE
%49S05,  %Variational principles of physics
%65Z05, % NUMERICAL ANALYSIS, Applications to physics
%74C05, %Small-strain, rate-independent theories (including rigid-plastic and elasto-plastic materials)
%74M15, % contact
%74R20. % anelastic fracture and damage
%}

\section{Introduction}\label{sect_Intro}
%        ~~~~~~~~~~~~

Adhesive contacts represent an important area in contact mechanics and
have numerous and continuously increasing applications.
The process of damaging the adhesive surfaces between bulk materials
is frequently referred to as {\it delamination} or debonding.
It is observed experimentally that sometimes, or rather typically,
more (or even substantially more) energy is needed in order for a
delamination to occur in the so-called Mode II (shear) than the respective
energy for a delamination in the so-called Mode I (opening). In general
when delamination proceeds in a mixed (and a-priori not known) mode,
we need a model which is sensitive to modes of delamination.

In this work, we use the Fr\'emond's concept analogous to a bulk damage,
assuming that the description of the damage is succeeded through a scalar
variable, which is defined along the adhesive interfaces,
taking values in the interval $[0,1]$
with 0 having the meaning of complete damage of the adhesive while 1
meaning complete operation of the adhesive, that is no damage appeared.
Moreover, we will consider the adhesive to have some elastic response,
also referred to as an imperfect or weak interface, opposite to the
rigid/ideally-brittle adhesive interface. Moreover, we will confine
ourselves to small strains and linearly responding
materials in the bulk. Actually, the present model, in its simplest form,
would correspond, following the classification in \cite{MarTru04IPFM},
to the so-called
``initially elastic Barenblatt model'' whose interface energy is given by a
convex quadratic function of displacement jump. It differs
form the classical Griffith model, which is not adequate for predicting
onset of delamination.

Let us remark that there is also another engineering model which
phenomenologically prescribes energy needed for delamination dependent on the
state (which is sometimes called ``non-associative'' model) through
the ratio of tangential and normal stresses or displacements
(whose {\sf arctan} is called a fracture-mode-mixity angle, see e.g.\
\cite{BanAsk00NFCI,HutSuo92MMCL,LieCha92ASIF,Mant08DRLM,SwLiLo99ITAM}).
Mathematical justification of such a model seems possible only
if a visco-elastic material with enough dissipative rheology (like
Kelvin-Voigt or Jeffrey) is considered, cf.\ \cite{RosRou13ACDM}.
Let us emphasize that this engineering model does not possess
any rigorous mathematical justification in
its typical engineering usage, i.e.\
the purely quasistatic elastic case when no viscosity is considered.

Microscopically, an additional dissipation in the Mode II
may be explained by a certain plastic process both in the adhesive
itself and in a narrow bulk vicinity
of the delamination surface before the actual delamination starts,
cf.~\cite{LieCha92ASIF,TveHut93IPMM}. Inspired by this,
another model has been devised in \cite{RoKrZe11DACM}
by introducing an extra inelastic parameter which describes some
{\it plastic-like interfacial slip} occurring possibly in the tangent
direction of an interface before its debonding is activated. This
interfacial plastification is not activated in Mode I, which allows
for dissipating some extra energy in Mode II compared to Mode I.
This model is associative in the sense that the
dissipation potential depends only on rates but not states.
Its rigorous analysis has been performed in \cite{RoMaPa13QMMD},
based on implicit discretisation in time and global minimization,
using the concept of globally stable energy-conserving (so-called energetic)
solutions devised by Mielke at al.\
\cite{Miel05ERIS,Miel11DEMF,MiRoSt08GLRR,MieThe04RIHM,MiThLe02VFRI}.

It is well known, however, that energetic solutions of rate-independent
problems governed by nonconvex energies (as inevitable in fracture
mechanics and in particular here too) tend to nonphysically too early jumps.
Instead of energy-driven and energy-conserving solution, some other
concepts seem more physically relevant, like vanishing-viscosity
solutions. See also the discussion about
energy versus stress or global versus local minimization
in mathematical literature \cite{CFMT00RBFE,MarTru04IPFM,Stef09VCRI}
and in engineering \cite{Legu02STCC,Mantic09ICOC}, and also the examples
\cite[Sect.\,9]{Cagn08VVAF} or \cite[Example\,7.1]{MiRoSa09MSJR}.
In general, all reasonable solutions fall into so-called
local solutions, invented \cite{ToaZan09AVAQ}, cf.\ also \cite{Miel11DEMF}.

Here, the undesired effect of too early jumps of globally stable 
energy-conserving solutions can be caused both by the influence of big energy
stored in the stressed bulk (cf.\ the explicit example in \cite{Roub13ACVE})
and also by a tendency to slide to
less dissipative mode of delamination (i.e.\ Mode I) even if the
direction of the traction stress would clearly tend to a more
dissipative mode (i.e. Mode II), as also indicated by numerical
experiments in \cite{RoMaPa13QMMD,RoKrZe11DACM}.

In the mode-mixity-insensitive model (i.e.\ Mode I dissipates
equally as Mode II), it have been observed in \cite{RoPaMA13QACV}
that the local solutions obtained by semi-implicit time discretisation
nicely coincides numerically with the vanishing-viscosity solutions
in all investigated examples; of course, the energy conservation has been
lost in such local but non-energetic solutions. Mathematical
justification of the semi-implicit time discretisation for
the quasistatic rate-independent problem was not scrutinized in
\cite{RoPaMA13QACV}, however.

The goal of this article is to devise a physically
relevant model (together with a corresponding
solution concept) for quasistatic mode-mixity-sensitive
delamination together with an efficient numerical strategy.
In Section~\ref{sect_AssocMod}, we briefly present the model
devised in \cite{RoKrZe11DACM,RoMaPa13QMMD} and in
Section~\ref{sect-LS} we define its solution using
the concept of local solutions from \cite{Miel11DEMF,ToaZan09AVAQ}.
Then, in  Section~\ref{sec-disc}, we devise a suitable semi-implicit
time discretisation and show its unconditional stability
in the sense that a-priori estimates can be proved, and then
prove convergence toward the local solutions of the
continuous problem. Eventually, in Section~\ref{sect-BEM},
we briefly present the fully discretised model
and outline its unconditional convergence if the time and space
discretisation refines, and present computational
simulations documenting this convergence as well as
physical relevancy of the model and its solution.
Let us emphasize that, in particular, it is for the first time
when the mode-mixity-sensitive delamination model and its solution
pursuing the 
maximum-dissipation principle and, due to Remark~\ref{rem-MDP} below,
the  stress-driven solution  concept 
in purely inviscid quasistatic situation
is analyzed as far as the convergence concerns
and, on a fixed discretisation, the approximate solutions can
efficiently be calculated non-iteratively at each time level by using a
linear-quadratic  programming algorithms.

\section{Quasistatic mode-mixity-sensitive delamination model}\label{sect_AssocMod}
%        ~~~~~~~~~~~~~~~~~~~~~~~~~~~~~~~~~~~~~~~~~~~~~~~

We will consider the evolution on a fixed finite time interval $[0,T]$
governed by a {\it stored energy} functional $\calE=\calE(t,u,z)$
and  a {\it dissipated energy} functional $\calR=\calR(\DT z)$
with the displacement field $u$ and an ``inelastic'' parameter field $z$
composed  here  from  delamination and interface-plasticity parameters.

The {\it delamination} (or {\it interfacial damage}) parameter is related to
fraction of adhesive bonds which are not broken. The
{\it interface-plasticity} parameter is motivated by the idea that,
microscopically,
the additional dissipation in Mode II may be explained by a certain plastic
processes both in the adhesive itself and in a narrow bulk vicinity of the
delaminating surface before the actual delamination starts \cite{TveHut93IPMM},
or by some rough structure of the interface \cite{ERDC90FEBI}.
In a certain idealization, these plastic processes are
more relevant in Mode II while do not manifest themselves significantly in
Mode I if the plastic strain is considered ‘incompressible’, i.e.\ trace-free.

Further,  we use the notation for the time
derivative $\DT z:=\frac{\mathrm d z}{\mathrm d t}$.
Specification of these energy functionals will be given later.
The general form of inclusions governing
the rate-independent evolution scrutinized in this article is the following
system of \emph{doubly nonlinear} degenerate abstract
\emph{static/evolution inclusions}, referred sometimes as Biot's equations
generalizing the original work \cite{Biot56TIT,Biot65MID}:
\begin{align}
\label{Gm}
\pl_u\calE(t,u,z)\ni0\quad\text{and}\quad
\pl\calR(\DT{z})+\pl_z\calE(t,u,z)\ni0,
\end{align}
where the symbol ``$\partial$'' refers to a (partial) subdifferential,
relying on that $\calR(\cdot)$, $\calE(t,\cdot,z)$, and
$\calE(t,u,\cdot)$ are convex functionals; the latter
inclusion in \eqref{Gm} thus contains the sum of two sets.

First we present in detail a plastic-type model with
kinematic-type hardening (like e.g.~in \cite{HanRed99PMTN,simo-hughes}) for
the above described delamination problem, devised,  analyzed,
and tested numerically in \cite{RoMaPa13QMMD,RoKrZe11DACM}.
The philosophy of the associative model is
to consider, besides some interface damage process described by a variable
$\zeta$, another inelastic process on the delaminating
surface $\GC$ which would be activated rather in fracture Mode II than in
Mode I, and thus more energy would be dissipated in Mode II than in
Mode I. This additional inelastic process involves an additional dissipative
variable $\pi$ having the meaning of the plastic-like
\emph{tangential slip} on $\GC$;  this variable defined on $\GC$
is a $(d{-}1)$-dimensional vector.  We will use a gradient theory
for   some  of the internal parameters used also, e.g.,\  in
\cite[Chap.14]{Frem02NST} or \cite{RoMaPa13QMMD,RoKrZe11DACM}.
In contrast to \cite{RoMaPa13QMMD,RoKrZe11DACM}, we consider here
the gradient of $\pi$ instead of $\zeta$ because now
we need strong convergence of all convex nonlinear terms,
which does not seem easy for a term like $|\nablaS\zeta|^r$ if
$\partial\calR$ is not bounded, as it is the case here because no re-bonding is
considered, i.e.\ only $\DT\zeta\le0$ is allowed,
cf.\ \eqref{R-delam-small-k-II}.

The relation to \eqref{Gm} is that $z=(\zeta,\pi)$ and, confining ourselves to
$\calR(\DT\zeta,\DT\pi)=\calR_0(\DT\zeta)+\calR_1(\DT\pi)$, the
system \eqref{Gm} takes the form
\begin{subequations}\label{Gm+}\begin{align}
\label{Gm1}
&\pl_u\calE(t,u,\zeta,\pi)\ni0,&&&&\text{(force equilibrium, Signorini contact)}
\\[-.3em]\label{Gm2}
&\pl\calR_0\big(\DT\zeta\big)+\pl_\zeta\calE(t,u,\zeta,\pi)
\ni0,&&&&\text{(a flow rule for interfacial damage)}
\\\label{Gm3}
&\pl\calR_1\big(\DT\pi\big)+\pl_\pi\calE(t,u,\zeta,\pi)
\ni0.&&&&\text{(a flow rule for interfacial plasticity)}
\end{align}\end{subequations}

To formulate the model, we consider two bounded Lipschitz domains
$\Omega_1,\Omega_2\subset\R^d$ ($d=2,3$)
with a common contact boundary $\GC=\partial\Omega_1\cap\partial\Omega_2$;
of course, the generalization for more than 2 bodies in contact is
straightforward.
Occasionally, we use the notation $\Omega=\Omega_1\cup\GC\cup\Omega_2$.
The contact boundary $\GC$ may undergo delamination.
We assume that the rest of the outer boundary $\partial\Omega$
is (up to $d{-}2$-dimensional zero-measure set)
the union of two disjoint open subsets $\GD$ and $\GN$
where the Dirichlet and the Neumann boundary conditions will be prescribed,
respectively. To ensure coercivity of the problem even after a possible 
complete delamination, we assume
\begin{align}\label{DQ-Dir}
\partial\Omega_1\cap\GD\ne\emptyset\qquad\&\qquad\partial\Omega_2\cap\GD\ne\emptyset.
\end{align}
On the Dirichlet part  $\GD$ of the boundary, we impose a time-dependent
displacement $\wD(t)$ and, on the boundary $\GN$, we impose a time-dependent 
traction $f(t)$.

The introduced associative delamination model is determined by the 
stored-energy functional
\begin{subequations}\label{delam-small-II}
\begin{align}
\label{E-delam-small-k-II}
& {\calE}(t,u,\zeta,\pi):=\!\left\{\begin{array}{ll}
   \displaystyle{\int_{\Omega\setminus\GC}\!\frac12\mathbb C
e(u){:}e(u)\,\d x-\!\int_{\GN}\!\!f(t){\cdot}u\,\d S\!\!\!}&\\[.5em]
\displaystyle{\ \ \ \
+\!\int_{\GC}\!\!
 \bigg(\zeta\Big(
\frac{\kappa_{_{\rm N}}}{2}\JUMP{u}{_{\rm N}}^{\!\!\!\!2}+
\frac{\kappa_{_{\rm T}}}{2}\big|\JUMP{u}{_{\rm T}}\!{-}\mathbb{T}\pi\big|^2\Big)
}\!\!\!&\\[.9em]
\displaystyle{\ \ \ \ \ \ \
+\,\frac{\kappa_{_{\rm H}}}{2}|\pi|^2
+\frac{\kappa_{_{\rm G}}}2\big|\nablaS\pi\big|^2 \bigg)}
\displaystyle{\,\d S}
&\!\!\!\!\!\!\!\!\!\!\!\!\!\!\!\text{if }u=\wD(t)\ \text{ on }\GD,\\[.0em]
&\!\!\!\!\!\!\!\!\!\!\!\!\!\!\!\JUMP{u}{_{\rm N}}\ge0,\
0\le\zeta\le1\text{ on }\GC,
\\[.3em]
   \infty&\!\!\!\!\!\!\!\!\!\!\!\!\!\!\!\text{elsewhere,}
\end{array}\right.\hspace*{-3em}
\intertext{and by the dissipated-energy functional}
\label{R-delam-small-k-II}
&\calR_0(\DT\zeta):=
\begin{cases}\displaystyle{\int_{\GC}\!\!a_{_{\rm I}}\big|\DT\zeta\big|\d S}
&\text{ if }\DT\zeta\le0\text{ a.e.~on }\GC,\\
\infty&\text{ otherwise},
\end{cases}
\qquad
\calR_1(\DT\pi):=\int_{\GC}\!\!\sigma_{\mathrm{yield}}\big|\DT\pi\big|\d S,
\end{align}
\end{subequations}
with $\mathbb C$ being the elastic-moduli tensor (possibly being $x$-dependent
and, in particular, may be different at the subdomains $\Omega_1$ and
$\Omega_2$), $e(u)=\frac12(\nabla u)^\top\!+\frac12\nabla u$ denoting the
small strain tensor, 
$a_{_{\rm I}}>0$ the prescribed phenomenological energy per unit area
dissipated (=\,needed for complete delamination,  often referred to as fracture energy or fracture toughness) in pure Mode I, 
$\nablaS$ a ``surface gradient'' (i.e.\ the tangential
derivative defined as $\nablaS v=\nabla v-(\nabla v{\cdot}\nu)\nu$ for
$v$ defined around $\GC$) and
$\JUMP{u}{}=\JUMP{u}{_{\rm N}}\nu+\JUMP{u}{_{\rm T}}$ with
$\JUMP{u}{_{\rm N}}=\JUMP{u}{}{\cdot}\nu$ with $\nu$ a unit normal to $\GC$;
in other words, $\nablaS v:=P(\nabla v)$ 
and $\JUMP{u}{_{\rm T}}=\JUMP{Pu}{}$ with the projector
$P={\mathbb I}-\nu\otimes\nu$ onto a tangent space. Here we
used the notation $\JUMP{u}{}$ for the differences of traces from both sides
of $\GC$. Note also that $\JUMP{u}{\rm N}$ is scalar valued while
$\JUMP{u}{_{\rm T}}$ is vector valued. Alternatively, pursuing the
concept of fields defined exclusively on $\Gamma$, we can consider
$v:\Gamma\to\R^d$ and extend it to a neighborhood of $\Gamma$ and
then again define $\nablaS v:=P(\nabla v)$ which, in fact, does
not depend on the particular extension. 
In fact, the model naturally does not depend on the chosen orientation.
 The value $\infty$ in \eqref{R-delam-small-k-II} guarantees that
$\DT\zeta\le0$ during the whole evolution everywhere on $\GC$, i.e.\
the interfacial damage evolution is irreversible, we can also say
uni-directional, in the sense that no re-bonding (i.e.\ no healing) is
allowed. In \eqref{E-delam-small-k-II}, 
$\mathbb{T}:\GC\to{\rm Lin}(\R^{d-1},\R^d)$ is to embed $\R^{d-1}$, where 
$\pi$ is valued, into $\R^d$, where $\JUMP{u}{_{\rm T}}$ has values,  so  that
$\JUMP{u}{_{\rm T}}\!{-}\mathbb{T}\pi$ has a sense,
 and we assume that $\mathbb{T}(x)(\R^{d-1})$ is 
the $(d{-}1)$-dimensional tangent space to $\GC$ at $x$  for
a.a.\ $x\in\GC$. The phenomenological elastic
constants $\kappa_{_{\rm N}}$ and  $\kappa_{_{\rm T}}$ in \eqref{E-delam-small-k-II}
describe the stiffnesses of linearly elastically responding adhesive in the
normal and tangential directions, respectively.
Typical phenomenology is that $\kappa_{_{\rm N}}$
is greater than $\kappa_{_{\rm T}}$; even, for isotropic adhesive,
a condition $\kappa_{_{\rm N}}/\kappa_{_{\rm T}}\ge2$ has been deduced in
\cite{TMGP11BEMA}, see also further references therein.

The unilateral constraint $\JUMP{u}{_{\rm N}}\ge0$ in \eqref{E-delam-small-k-II}
guarantees  infinitesimal 
nonpenetration  before and   after delamination (the so-called Signorini contact)
and impossibility of delamination by  pure  compression.
To produce desired effects, the model should work with parameters
satisfying\\[-1.7em]
\begin{align}\label{2-sided-condition}
\ \ \frac12\kappa_{_{\rm T}}a_{_{\rm I}}<\sigma_{\mathrm{yield}}^2\le2\kappa_{_{\rm T}}a_{_{\rm I}}.
\end{align}
{M}ore specifically, the upper bound of the yield stress is necessary for
making possible to initiate plastic slip before the total interface damage,
while the lower bound is required to avoid plastic slip
evolution after complete debonding when $\zeta=0$. 
Then, one can see that the overall dissipated energy in Mode II, denoted by
$a_{_{\rm II}}$, is\\[-1.7em]
\begin{align}\label{ch5:delam-small-a-II}
\ \ a_{_{\rm II}}=a_{_{\rm I}}+
\frac{2\kappa_{_{\rm T}}a_{_{\rm I}}-\sigma_{\mathrm{yield}}^2}{ 2 \kappa_{_{\rm H}}}
\end{align}
cf.\ \cite{RoMaPa13QMMD,RoKrZe11DACM}  where, however, the contribution of
the hardening after complete delamination
was ignored. This hardening energy, although being a part of the stored
energy, cannot be gained back  (assuming (2.5))  and thus is effectively dissipated for ever after
the delamination in Mode II is completed. 
For example, for
$\sigma_{\mathrm{yield}}=\sqrt{\kappa_{_{\rm T}}a_{_{\rm I}}}$,
\eqref{2-sided-condition} is satisfied and
$a_{_{\rm II}}/a_{_{\rm I}}=1+\kappa_{_{\rm T}}/(2\kappa_{_{\rm H}})$.
In particular, by choosing
$\kappa_{_{\rm H}}>0$ small, this model can handle arbitrarily large
ratio $a_{_{\rm II}}/a_{_{\rm I}}$; let us emphasize that in engineering
situations, this ratio is often over 10.

The typical occurrence of jumps of solutions needs a careful definition
relying on the time derivative $\calE_t'(\cdot,u,\zeta,\pi)$ for 
$(u,\zeta,\pi)$ fixed, cf.\ the last term in \eqref{def-ls-mix-engr}. This 
obviously requires $\wD$ in \eqref{E-delam-small-k-II} constant in time
 to avoid the situation that, for $u$ fixed such that
$u|_{\GD}=\wD(t_0)$, the value $\calE(t,u,\zeta,\pi)$ is finite for $t=t_0$ while
it equals $\infty$ for $t\ne t_0$ and thus $\calE_t'(\cdot,u,\zeta,\pi)$
cannot exist at $t=t_0$. In a general case  $\DTwD\ne0$,
we make a substitution of $u{+}\uD(t)$ with
$\uD(t)$ being a suitable extension of $\wD(t)$. Then, up to the time-dependent
constant $\int_{\Omega\setminus\GC}\frac12\mathbb Ce(\uD(t)){:}e(\uD(t))
\,\d x-\int_{\GN}f(t){\cdot}\uD(t)\,\d S$, \eqref{E-delam-small-k-II} is to be
replaced by
\begin{subequations}\label{E-delam-small-k-II+}\begin{align}
&\!\!\!\!{\calE}(t,u,\zeta,\pi):=\!\left\{\begin{array}{ll}
   \displaystyle{\int_{\Omega\setminus\GC}\!
\frac12\mathbb Ce(u){:}e(u)\,\d x-\big\langle f_1(t),u\rangle
\hspace{-0em}}&\\[-.4em]\displaystyle{\ \ \ \ \ \
+\!\int_{\GC}\!\!
 \bigg(\zeta\Big(
\frac{\kappa_{_{\rm N}}}{2}\JUMP{u}{_{\rm N}}^{\!\!\!\!2}+
\frac{\kappa_{_{\rm T}}}{2}\big|\JUMP{u}{_{\rm T}}\!{-}\mathbb{T}\pi\big|^2\Big)
}\hspace{-2em}&\\[-.0em]
\displaystyle{\ \ \ \ \
+\,\frac{\kappa_{_{\rm H}}}{2}|\pi|^2
+\frac{\kappa_{_{\rm G}}}2\big|\nablaS\pi\big|^2 \bigg)}
\displaystyle{\,\d S}
&\!\!\!\!\!\!\!\!\text{if }u=0\text{ on }\GD\,\text{ and }\,
\\[.0em]&\!\!\!\!\!\!\!\!
\JUMP{u}{_{\rm N}}\ge0,\ 0\le\zeta\le1\text{ on }\GC,
\\[.3em]
   \infty&\!\!\!\!\!\!\!\!\text{elsewhere,}
\end{array}\right.\hspace*{-3em}
\\\label{E-delam-f1}
&\qquad\qquad\qquad\quad\text{ with }\ \big\langle f_1(t),v\rangle=
\int_{\GN}\!\!f(t){\cdot}v\,\d S
-\int_{\Omega\setminus\GC}\!\mathbb Ce(\uD(t)){:}e(v)\,\d x.
\end{align}\end{subequations}
In fact, we have assumed that $\GD$ is far from $\GC$ so that we can
have $\uD(t)|_{\GC}=0$ not to affect the integral over $\GC$ in
\eqref{E-delam-small-k-II} by the shift $u\mapsto u{+}\uD(t)$.

An alternative way is to avoid this transformation by considering
only the trace of $u$ on $\GC$ as the state variable instead of $u$. This is
possible by using the boundary-integral-equation  (BIE) 
method  which
evaluates the bulk integral and eliminates the constraint $u|_{\GD}=\wD(t_0)$
in \eqref{E-delam-small-k-II} by solving the boundary-value problem
governed by minimization of this integral under the
condition that $u$ is prescribed on $\GD\cup\GC$.
After spatial discretisation, BIE becomes the boundary-element method,
which is in fact how we implement the problem
in Section~\ref{sect-BEM} below, although the analysis
is performed on the more conventional base of the
transformed functional \eqref{E-delam-small-k-II+}.

We will consider an initial-value problem for the system \eqref{Gm+} by
prescribing
\begin{align}\label{IC}
u(0)=u_0,\qquad\zeta(0)=\zeta_0,\qquad\pi(0)=\pi_0.
\end{align}

\section{Local solutions}\label{sect-LS}
%        ~~~~~~~~~~~~~~~

We will also abbreviate the time interval $I=[0,T]$ with $T>0$ a fixed
time horizon, and $\SC=I\times\GC$.

We will use the standard notation $W^{1,p}(\Omega)$ for the Sobolev space
of functions having the gradient in the Lebesgue space $L^p(\Omega;\R^d)$. If valued in
$\R^n$ with $n\ge2$, we will write $W^{1,p}(\Omega;\R^n)$, and furthermore we
use the shorthand notation $H^1(\Omega;\R^n)=W^{1,2}(\Omega;\R^n)$.
Similarly, we will use Lebesgue and Sobolev space on the
$(d{-}1)$-dimensional manifold $\GC$, assumed Lipschitz so that 
a local rectification for defining the surface gradient
$\nablaS$  can be performed a.e.\ on $\GC$.
We also use the notation of ``$\,\cdot\,$'' and ``$\,:\,$''
for a scalar product of vectors and 2nd-order tensors, respectively.
For a Banach space $X$, $L^p(I;X)$ will denote the Bochner space of
$X$-valued Bochner measurable functions $u:I\to X$ with its norm
$\|u(\cdot)\|$ in $L^p(I)$, here $\|\cdot\|$ stands for
the norm in $X$. Further, $W^{1,p}(I;X)$ denotes the Banach space of mappings
$u:I\to X$
whose distributional time derivative is in $L^p(I;X)$, while
${\rm BV}(I;X)$ will denote the space of mappings $u:I\to X$
with a bounded variations, i.e.\ $
\sup_{0\le t_0<t_1<...<t_{n-1}<t_n\le T}\sum_{i=1}^n\|u(t_i){-}u(t_{i-1})\|<\infty$
where the supremum is taken over all finite partitions of the interval
$I=[0,T]$. By  ${\rm B}(I;X)$ we denote the space of bounded measurable
(everywhere defined) mapping $I\to X$.

The concept of local solutions has been introduced for a special crack
problem in \cite{ToaZan09AVAQ}  and independently also in
\cite{Stef09VCRI},  and further generally investigated in
\cite{Miel11DEMF}. Here, we additionally combine it with the concept of
semi-stability as invented in \cite{Roub09RIPV}.
We adapt the general definition directly to our specific problem,
which will lead to two semi-stability conditions for $\zeta$ and
$\pi$, respectively:

\begin{definition}[Local solutions]\label{def-LS}
We call a measurable mapping
$(u,\zeta,\pi):I\to H^1(\Omega{\setminus}\GC;\R^d)\times
L^\infty(\GC)\times H^1(\GC;\R^{d-1})$ a local solution to the
delamination problem \eqref{Gm+}--\eqref{IC} if the initial conditions
\eqref{IC} are satisfied, if $\JUMP{u}{_{\rm N}}\ge0$ on
$\SC$ and, for some $J\subset I$ at most countable (containing
time instances where the solution may possibly jump), it holds that:
\begin{subequations}\label{def-mix-LS}
\begin{align}
\nonumber
&\!\!\!\forall t\!\in\!I{\setminus}J\ \forall v\!\in\!
 H^1(\Omega{\setminus}\GC;\R^d),\ \JUMP{v}{_{\rm N}}\ge0:\ \ \
\\\label{def-mix-ls-VI}
&\qquad\qquad\qquad
\big\langle\pl_u\calE\big(t,u(t),\zeta(t),\pi(t)\big),v{-}u(t)\big\rangle\ge0,
\\\nonumber
&\!\!\!\forall t\!\in\!I\ \forall\widetilde\zeta\!\in\! L^\infty(\GC),\ \,
0\!\le\!\widetilde\zeta\!\le\!\zeta(t)\text{ a.e.\ on }\GC:\ \
\\\label{def-ls-mix-semi-stab2}
&\qquad\qquad\qquad
\calE\big(t,u(t),\zeta(t),\pi(t)\big)\le\calE\big(t,u(t),\widetilde\zeta,\pi(t)\big)
+\calR_0\big(\widetilde\zeta{-}\zeta(t)\big),
\\\nonumber
&\!\!\!\forall t\!\in\!I{\setminus}J\ \forall\widetilde\pi\!\in\! H^1(\GC;\R^{d-1}):\quad
\\&\qquad\qquad\qquad
\calE\big(t,u(t),\zeta(t),\pi(t)\big)\le\calE\big(t,u(t),\zeta(t),\widetilde\pi\big)
+\calR_1\big(\widetilde\pi{-}\pi(t)\big),
\label{def-ls-mix-semi-stab1}
\\\nonumber&\!\!\!\forall0\!\le\!t_1\!\le\!t_2\!\le\!T:\ \ \
\calE\big(t_2,u(t_2),\zeta(t_2),\pi(t_2)\big)
+\Diss_{\calR_1}^{}\big(\pi;[t_1,t_2]\big)
\\[-.5em]\label{def-ls-mix-engr}&\hspace{3em}
+
\int_{\GC}\!\!\!a_{_{\rm I}}\big(\zeta(t_1){-}\zeta(t_2)\big)\,\d S
\le\calE\big(t_1,u(t_1),\zeta(t_1),\pi(t_1)\big)-\!\int_{t_1}^{t_2}\!\!
\big\langle\DT f_1,u\big\rangle\,\d t
\end{align}\end{subequations}
where $f_1$ is from \eqref{E-delam-f1} and
$\Diss_{\calR_1}^{}(\pi;[r,s]):=\sup\sum_{j=1}^N
\sigma_{\mathrm{yield}}|\pi(t_{j-1}){-}\pi(t_j)|$ with the supremum taken over all 
finite partitions 
$r\leq t_0^{}\!<\!t_1^{}\!<\!\cdots\!<\!t_{N-1}^{}\!<\!t_N^{}\!\leq s$.
\end{definition}

Let us comment the above definition briefly. Obviously, \eqref{Gm1} means
precisely \eqref{def-mix-ls-VI},  which more in detail here means 
that
$\int_{\Omega{\setminus}\GC}\mathbb{C} e(u(t)){:}e(v-u(t))\,\d x
+\int_{\GC}\zeta(t)(\kappa_{_{\rm N}}\JUMP{u(t)}{_{\rm N}},
\kappa_{_{\rm T}}\JUMP{u(t)}{_{\rm T}}\!-\mathbb{T}\pi)
{\cdot}\JUMP{v-u(t)}{}\d S\ge\langle f_1(t),v-u(t)\rangle$
for all $v\!\in\!H^1(\Omega{\setminus}\GC;\R^d)$
with $\JUMP{v}{_{\rm N}}\ge0$.
 Note that \eqref{def-mix-ls-VI} specifies also the boundary conditions
for $u$, namely $u=0$ on $\GD$ because otherwise
$\calE( t, u,\zeta,\pi)=\infty$ would violate \eqref{def-mix-ls-VI}
for $v$ which satisfies $v=0$ on $\GD$, and also
$\nu^\top\mathbb C e(u)=f$ on $\GN$ can be proved by standard
arguments based on Green's theorem. 
As $\calR_1$ is homogeneous degree-1, always
$\partial\calR_1(\DT\pi)\subset\partial\calR_1(0)$ and thus \eqref{Gm3} implies
$\partial\calR_1(0)+\partial_\pi\calE(u,\zeta,\pi)\ni0$.
From the convexity of $\calR_1$ when taking into account
that $\calR_1(0)=0$, the latter inclusion is equivalent to
$\calR_1(v)+\langle\partial_\pi\calE(u(t),\zeta(t),\pi(t)),v\rangle\ge0$
for any $v\in H^1(\GC;\R^{d-1})$. Substituting $v=\widetilde z-z(t)$ and using
the convexity of
$\calE(t,u,\zeta,\cdot)$, we obtain the \emph{semi-stability}
 \eqref{def-ls-mix-semi-stab1}  of $\pi$ at time $t$.
Analogously, we obtain also \eqref{def-ls-mix-semi-stab2} from  \eqref{Gm2}.
Eventually, \eqref{def-ls-mix-engr} is the (im)balance of the mechanical
energy with the last term representing a ``complementry'' work of
external forces arising from the usual work by a by-part integration in time. 
This generalizes the standard definition of the {\it weak solution}
to \eqref{Gm+} to the case when $\calE(t,\cdot,\cdot,\cdot)$ is not
smooth, cf.\ \cite{Roub??MDLS} for details.

To be more precise, the concept of local solutions as used in
\cite{Miel11DEMF,ToaZan09AVAQ} requires $J$ only to have a zero Lebesgue
measure and also \eqref{def-ls-mix-semi-stab2} is valid
only for a.a.\ $t$. On the other hand, conventional weak solutions
allow even \eqref{def-ls-mix-engr} holding only for a.a.\ $t_1$ and $t_2$.
Later, our approximation method will provide convergence to this
slightly stronger local solutions, which motivates us to have tailored
Definition~\ref{def-LS} straight to our results.

Actually, local solutions form essentially the largest
reasonable class of solutions for \eqref{Gm}, coinciding
(in the above mentioned weaker form) with the conventional weak solutions,
cf.\ \cite{Roub??MDLS}. It includes
the mentioned energetic solutions \cite{Miel05ERIS,MieThe04RIHM},
the vanishing-viscosity solutions, the balanced-viscosity (so-called BV)
solutions, parametrized solutions, etc.; cf.\ \cite{Miel11DEMF,MieRou15RIST} 
for a survey, and also stress-driven-like 
solutions obeying maximum-dissipation principle
 in some sense. The energetic solution has often tendency
to rupture unphysically early and rather in the less dissipative
Mode I even if there should be rather Mode II expected; cf.\
\cite{VoMaRo??BEIM} for a comparison on several computational experiments.
The approximation method we will use in this article leads rather to 
the stress-driven  option, cf.\ Remarks~\ref{rem-MDP} and \ref{rem-DMP} below.

Anyhow, let us mention that, in \cite{RoMaPa13QMMD},
existence of the globally stable energy-conserving local
solutions of this model has been proved under the following assumptions:
\begin{subequations}\label{ass}\begin{align}\label{ass-C}
&\mathbb C^{(i)}\text{ positive definite, symmetric},\ \
\kappa_{_{\rm G}},\kappa_{_{\rm H}}>0,\ \ \kappa_{_{\rm T}},\kappa_{_{\rm N}}\ge0,\ \
a_{_{\rm I}},\sigma_{\mathrm{yield}}>0,
\\\label{ass-BC}
&\wD\in W^{1,1}(0,T;W^{1/2,2}(\GD;\R^d)),\ \ \
\\[-.3em]\label{ass-BC+}
&f\in W^{1,1}(0,T;L^p(\GN;\R^d))\ \ \ \text{ with }
\ \ p\begin{cases}>1&\text{for }d=2,\\[-.2em]
=2{-}2/d\!&\text{for }d\ge3\end{cases}
\\\label{ass-IC}
&(u_0,\zeta_0,\pi_0)\in H^1(\Omega{\setminus}\GC;\R^d){\times}
L^\infty(\GC){\times}H^1(\GC;\R^{d-1}),
\\
&\forall (\widetilde u,\widetilde\zeta,\widetilde\pi):
\ \ \ \calE(0,u_0,\zeta_0,\pi_0)\le
\calE(0,\widetilde u,\widetilde\zeta,\widetilde\pi)
+\calR(\widetilde\zeta{-}\zeta,\widetilde\pi{-}\pi).
\end{align}\end{subequations}
{T}he last condition, called stability at $t=0$, is
needed to ensure energy conservation and
will not be needed for general local solutions. 
The qualification \eqref{ass-BC} allows for an extension
$\uD$ of $\wD$ which belongs to $W^{1,1}(0,T;H^1(\Omega;\R^d))$;
in what follows, we will consider some extension with this property.

\begin{remark}[{M}aximum-dissipation principle]
\label{rem-MDP}
\upshape
The degree-1 homogeneity of $\calR_0$ and $\calR_1$ defined in
\eqref{R-delam-small-k-II} allows for further interpretation of the flow
rules \eqref{Gm2} and \eqref{Gm3}. Using maximal-monotonicity of the
subdifferential, \eqref{Gm3} means just that
$\langle\widetilde{\mathfrak{f}}-\mathfrak{f},v-\DT\pi\rangle\ge0$
for any $v$ and any $\widetilde{\mathfrak{f}}\in\partial\calR_1(v)$
with the  available  driving force
$\mathfrak{f}\in-\partial_\pi\calE(t,u,\zeta,\pi)$; the
adjective ``available'' becomes sensible especially if
$\partial_\pi\calE(t,u,\zeta,\pi)$ is set-valued
because not all available $\mathfrak{f}$'s are compatible with
$\mathfrak{f}\!\in\!\partial\calR_1(\DT\pi)$ and
can be realized during evolution.
In particular, for $v=0$, defining the convex set
$K_1:=\partial\calR_1(0)$, one obtains
\begin{subequations}\label{e1:max-dis-principle}
\begin{align}\label{e1:max-dis-principle-a}
\big\langle \mathfrak{f}(t),\DT\pi(t)\big\rangle
=\max_{\widetilde{\mathfrak{f}}\in K_1}\big\langle\widetilde{\mathfrak{f}},
\DT\pi(t)\big\rangle\ \ \ \text{ with some }\ \ \mathfrak{f}(t)\in-\partial_\pi\calE(t,u(t),\zeta(t),\pi(t)).
\end{align}
To derive it, we have used that
$\mathfrak{f}\in\partial\calR_1(\DT\pi)\subset\partial\calR_1(0)=K_1$
thanks to the degree-0 homogeneity of
$\partial\calR_1$, so that always $\langle \mathfrak{f},\DT\pi\rangle
\le\max_{\widetilde{\mathfrak{f}}\in K_1}\langle\widetilde{\mathfrak{f}},\DT\pi\rangle$. The identity \eqref{e1:max-dis-principle-a} says that the
dissipation due to the driving force $\mathfrak{f}$
is maximal provided that the order-parameter rate $\DT\zeta$ is 
kept fixed, while the vector of possible driving
forces $\widetilde{\mathfrak{f}}$ varies freely over all admissible driving
force from $K_1$. This just resembles the so-called Hill's
\emph{maximum-dissipation principle}  articulated just for plasticity in
 \cite{Hill48VPMP}. Also it says that the rates are orthogonal
to the ``elastic domain'' $K_1$, known as an
orthogonality principle \cite{Zieg58AGOP} generalizing Onsager's
principle \cite{Onsa31RRIP}. See also
\cite{HacFis07RPMD,Lub84MDPG,RajSri98MIBM,ZieWeh87PMER}.
 Actually, R.\,Hill \cite{Hill48VPMP} used it for a situation
where $\calE(t,\cdot)$ is convex while, in a general nonconvex case
as also here, it holds only along absolutely continuous paths
(i.e.\ in stick or slip regimes) which are sufficiently
regular in the sense $\DT\pi$ is valued not only in 
$L^1(\GC;\R^{d-1})$ but also in $H^1(\GC;\R^{d-1})^*$ 
but it does  certainly  not need to hold during jumps. 
Analogously it holds also for $\zeta$, defining
$K_0:=\partial\calR_0(0)$, i.e.
\begin{align}\label{e1:max-dis-principle-b}
\big\langle \mathfrak{g}(t),\DT\zeta(t)\big\rangle
=\max_{\widetilde{\mathfrak{g}}\in K_0}\big\langle\widetilde{\mathfrak{g}},
\DT\zeta(t)\big\rangle\ \ \ \text{ with some }\ \ 
\mathfrak{g}(t)\in-\partial_\zeta\calE(t,u(t),\zeta(t),\pi(t)).
\end{align}\end{subequations}
As $\calE(t,u,\zeta,\cdot)$ is smooth,
the maximum-dissipation relation \eqref{e1:max-dis-principle-a}
written in the form
\begin{align*}
\langle-\calE'_\pi(t,u(t),\zeta(t),\pi(t)), \DT \pi(t)\rangle=\max\langle K_1, \DT \pi(t)\rangle=\calR_1(\DT\pi(t))
\end{align*}
summed with the semistability \eqref{def-ls-mix-semi-stab1}
which can be written in the form
\begin{align*}
\calR_1(\widetilde\pi)
+\langle\calE'_\pi(t,u(t),\zeta(t),\pi(t)),\widetilde\pi\rangle\ge0
\end{align*}
thanks to the convexity of $\calE(t,u,\zeta,\cdot)$ yields
\begin{align}
\calR_1(\widetilde\pi)
+\langle\calE'_\pi(t,u(t),\zeta(t),\pi(t)),\widetilde\pi- \DT \pi(t) \rangle
\ge\calR_1(\DT\pi(t))
\end{align}
for any $\widetilde\pi$, which just means that
$\mathfrak{f}(t)=-\calE'_\pi(t,u(t),\zeta(t),\pi(t))
\in\partial\calR_1(\DT\pi(t))$. This exactly means that the evolution
of $\pi$ is {\it governed by a thermodynamical driving force}
$\mathfrak{f}$  (we say that it is ``stress-driven'')  and it reveals
the role of the  maximum-dissipation principle in combination
with semistability. Using the convexity of $\calE(t,u,\cdot,\pi)$,
a similar argument can be applied for
\eqref{e1:max-dis-principle-b} in combination with semistability
\eqref{def-ls-mix-semi-stab2} even if $\calE(t,u,\cdot,\pi)$ is not
smooth. 
\end{remark}

\begin{remark}[Integrated maximum-dissipation principle]
\label{rem-IMDP}
\upshape
Let us emphasize that, in general, $\DT\zeta$ and $\DT\pi$ are measures
possibly having singular parts concentrated at rupture times where the
solution and also the driving forces need not be continuous.
Even if  $\DT\zeta$ and $\DT\pi$ are absolutely continuous, in
our infinite-dimensional case the driving forces need not be in duality
 with them, as already mentioned in Remark~\ref{rem-MDP}. 
So \eqref{e1:max-dis-principle} is analytically not justified in any sense. 
For this reason,
an {\it Integrated} version of the {\it Maximum-Dissipation
Principle} (IMDP) was devised in \cite{Roub??MDLS} for a bit simpler
case involving only one maximum-dissipation relation.
Realizing that
$\max_{\widetilde{\mathfrak{f}}\in K_1}\langle\widetilde{\mathfrak{f}},
\DT\pi\rangle=\calR_1(\DT\pi)$ and
similarly $\max_{\widetilde{\mathfrak{g}}\in K_0}\langle\widetilde{\mathfrak{g}},
\DT\zeta\rangle=\calR_0(\DT\zeta)$, the integrated version of
\eqref{e1:max-dis-principle} reads here as:
\begin{subequations}\label{IMDP}\begin{align}\label{IMDP-a}
\!\int_{t_1}^{t_2}\!\!\!\mathfrak f(t)\,\d  \pi (t)
=\!\int_{t_1}^{t_2}\!\!\!\calR_1(\DT\pi)\,\d t\ \ \ 
\ \ \text{ with some }\
\mathfrak f(t)\!\in\!-\partial_\pi\calE(t,u(t),\zeta(t),\pi(t)),\
\\\label{IMDP-b}
\!\!\int_{t_1}^{t_2}\!\!\!\mathfrak g(t)\,\d  \zeta  (t)
=\!\int_{t_1}^{t_2}\!\!\!\calR_0(\DT\zeta)\d t\ \ \ 
\ \ \text{ with some }\ \:
\mathfrak g(t)\!\in\!-\partial_\zeta\calE(t,u(t),\zeta(t),\pi(t))\
\end{align}\end{subequations}
to be valid for any $0\le t_1<t_2\le T$. This definition is inevitably
a bit technical and, without sliding into too much details, let us
only mention that the left-hand-side integrals in \eqref{IMDP}
are the so-called {\it lower Riemann-Stieltjes integrals} defined by
suprema of lower Darboux sums, i.e.\ in the case \eqref{IMDP-a} as
\begin{align*}
\int_r^s\!\mathfrak f(t)\,\d  \pi (t):=\!\!
\sup_{\substack{N\in\N\\r=t_0^{}<t_1^{}<...<t_{N-1}^{}<t_N^{}=s^{^{}}}}\,
\sum_{j=1}^{N}\,\inf_{t\in[t_{j-1},t_j]}
\big\langle\mathfrak f(t), \pi (t_j){-} \pi (t_{j-1})\big\rangle,
\end{align*}
while the right-hand-side integrals are just the integrals of measures and
equal to $\Diss_{\calR_1}(\pi;[t_1,t_2])$ and $\Diss_{\calR_0}(\zeta;[t_1,t_2])$,
respectively. The IMDP \eqref{IMDP} is satisfied on any interval $[t_1,t_2]$
where the solution is absolutely continuous  with sufficiently regular time 
derivatives; then the integrals in \eqref{IMDP}
are the conventional Lebesgue integrals, in particular the left-hand
sides in \eqref{IMDP} are  
$\int_{t_1}^{t_2}\langle\mathfrak f(t),\DT\pi(t)\rangle\,\d t$
and $\int_{t_1}^{t_2}\langle\mathfrak g(t),\DT\zeta(t)\rangle\,\d t$, 
respectively. 
The particular importance of IMDP is  especially  at jumps, i.e.\ at 
times when abrupt delamination possibly happens.
It is shown in \cite{MieRou15RIST,Roub??MDLS} on various
finite-dimensional examples of ``damageable springs'' that this
IMDP can identify too early rupturing local solutions  when the
driving force is obviously unphysically low (which occurs quite typically  
in particular  within  the energetic  solutions of systems
governed by nonconvex potentials like here) and its satisfaction for 
left-continuous local solutions indicates that the evolution is stress driven,
as explained in Remark~\ref{rem-MDP}.
On the other hand, it does not need to be satisfied
even in physically well justified stress-driven local solutions. For example,
it happens if two springs with different fracture toughness organized
in parallel rupture at the same time (although even in this situation
our algorithm \eqref{delam-small-II-LS-min}  below will  give 
a correct approximate solution).
Therefore,  even the IMDP \eqref{IMDP} may serve only as a 
sufficient aposteriori condition whose satisfaction verifies 
the obtained local solution as a physically relevant in the sense
that it is stress driven but its dissatisfaction does not mean anything.
Moreover, 
we will rely rather on some approximation of IMDP, as described in
Remark~\ref{rem-DMP} below.
\end{remark}

\section{Semi-implicit time discretisation, its stability and convergence}
%        ~~~~~~~~~~~~~~~~~~~~~~~~~~~~~~~~~~~~~~~~~~~~~~~~~~~
\label{sec-disc}

To prove existence of the physically relevant solution, we use
a constructive method relying on time discretisation and
the weak compactness of level sets of the minimization problems
arising at each time level. When further discretised in space, it will later
in Sect.~\ref{sect-BEM} yield a computer implementable efficient
algorithm.

For the mentioned time discretisation, we use an equidistant partition
of the time interval $I=[0,T]$ with a time step $\tau>0$, assuming
$T/\tau\in\N$, and denote $\{u_\tau^k\}_{k=0}^{T/\tau}$ an approximation
of the desired values $u(k\tau)$, and similarly $\zeta_\tau^k$ is to
approximate $\zeta(k\tau)$, etc.

We use a decoupled  semi-implicit time discretisation  with the
fractional steps  based on the splitting
of the state variables governed by the separately-convex
character of $\calE(t,\cdot,\cdot,\cdot)$. This will make the
numerics considerably easier than any other splitting and simultaneously
 may lead  to a physically relevant solutions governed rather by
stresses (if the maximum-dissipation principle  holds at least approximately
in the sense of Remark~\ref{rem-DMP} below) than by energies and will prevent 
too-early debonding, as already announced in Section \ref{sect_Intro}. 
More specifically, exploiting the convexity of both
$\calE(t,\cdot,\zeta,\cdot)$ and $\calE(t,u,\cdot,\pi)$,
this splitting will be considered as $(u,\pi)$
and $\zeta$. This yields alternating convex minimization.
Thus, for $(\zeta_\tau^{k-1},\pi_\tau^{k-1})$ given, we obtain two
minimization problems
\begin{subequations}\label{delam-small-II-LS-min}
\begin{align}\label{minimize-1}
\left.\begin{array}{ll}
\text{minimize}&\calE(k\tau,u,\zeta_\tau^{k-1},\pi)+\calR_1(\pi{-}\pi_\tau^{k-1})
\\\text{subject to}&(u,\pi)\in H^1(\Omega{\setminus}\GC;\R^d)\times H^1(\GC;\R^{d-1}),
\end{array}\ \right\}
\intertext{and, denoting the unique solution as $(u_\tau^k,\pi_\tau^k)$,}
\label{minimize-2}
\left.\begin{array}{ll}
\text{minimize}&\calE(k\tau,u_\tau^k,\zeta,\pi_\tau^k)
+\calR_0(\zeta{-}\zeta_\tau^{k-1})
\\\text{subject to}&\zeta\in L^\infty(\GC),\ \ 0\le\zeta\le1,
\end{array}\hspace{4em}\right\}
\end{align}
\end{subequations}
and denote its (possibly not unique) solution by $\zeta_\tau^k$.

Existence of the discrete solutions $(u_\tau^k,\zeta_\tau^k,\pi_\tau^k)$ is 
straightforward by the mentioned compactness arguments. Rather, it is important
that both problems \eqref{delam-small-II-LS-min} have the linear-quadratic
structure, the former one after applying the Mosco-type transformation,
cf.\ \cite[Lemma~4]{Roub02EMMP}. This obviously facilitates their numerical 
treatment; cf.\ Section~\ref{sect-BEM} below.

We define the piecewise-constant interpolants
\begin{align}\label{e3:def-of-interpolants}
\left.\begin{array}{l}
\bar u_\tau(t)= u_\tau^k\,\ \ \ \&\ \ \ \underline u_\tau(t)= u_\tau^{k-1},
\\[.2em]\bar\zeta_\tau(t)=\zeta_\tau^k\,\ \ \ \&\ \ \
\underline\zeta_\tau(t)=\zeta_\tau^{k-1},
\\[.2em]\bar\pi_\tau(t)=\pi_\tau^k\ \ \ \&\ \ \ \underline\pi_\tau(t)=\pi_\tau^{k-1},
\\[.2em]\bar\calE_\tau(t,u,\zeta,\pi)=\calE(k\tau,u,\zeta,\pi)
\end{array}\ \right\}\text{ for }(k{-}1)\tau<t\le k\tau.
\end{align}
Later in Remark~\ref{rem-DMP}, we will use also the piecewise affine
interpolants
\begin{align}\label{e3:def-of-interpolants+}
\left.\begin{array}{l}
\zeta_\tau(t)=\frac{t-(k{-}1)\tau}\tau\zeta_\tau^k
+\frac{k\tau-t}\tau\zeta_\tau^{k-1},
\\[.2em]\pi_\tau(t)=\frac{t-(k{-}1)\tau}\tau\pi_\tau^k
+\frac{k\tau-t}\tau\pi_\tau^{k-1}
\end{array}\ \right\}\text{ for }(k{-}1)\tau<t\le k\tau.
\end{align}

The important attribute of the discretisation \eqref{delam-small-II-LS-min}
is also its numerical stability and satisfaction of a suitable discrete analog
of \eqref{def-mix-LS}, namely:

\begin{proposition}[Stability of the time discretisation]
Let (\ref{ass}a-d) hold and, in terms of the interpolants 
\eqref{e3:def-of-interpolants}, $(\bar u_\tau,\bar\zeta_\tau,\bar\pi_\tau)$ be 
an approximate solution obtained by \eqref{delam-small-II-LS-min}.
Then, the following a-priori estimates holds
\begin{subequations}\label{e6:delam-mix-ls-estimates}
\begin{align}\label{e6:delam-mix-ls-estimate1}
&\big\|\bar u_\tau\big\|_{ L^\infty(I; H^1(\Omega{\setminus}\GC;\R^d))}
\le C,&&
\\\label{e6:delam-mix-ls-estimate2}
&\big\|\bar\zeta_\tau\big\|_{ L^\infty(\SC)
\cap{\rm BV}(I; L^1(\GC))}\le C,
\\\label{e6:delam-mix-ls-estimate3}
&\big\|\bar\pi_\tau\big\|_{ L^\infty(I; H^1(\GC;\R^{d-1}))
\cap{\rm BV}(I; L^1(\GC;\R^{d-1}))}\le C.
&&
\end{align}\end{subequations}
Moreover, the obtained approximate solution satisfies
for any $t\in I$ the variational inequality for the displacement:
\begin{subequations}\label{e6:delam-mix-ls}
\begin{align}\nonumber
&\forall\widetilde u\in H^1(\Omega{\setminus}\GC;\R^d),\ \JUMP{\widetilde u}{_{\rm N}}\!\ge0:
\\&\qquad\qquad\quad
\big\langle\pl_u\calE\big(t,\bar u_\tau(t),\underline\zeta_\tau(t),\bar\pi_\tau(t)\big), \widetilde u{-}\bar u_\tau(t)\big\rangle\ge0,
\label{e6:delam-mix-ls-VI}
\intertext{with $t_\tau\!:=\min\{k\tau\!\ge\!t;\ k\!\in\!\N\}$,
two separate semi-stability conditions for $\bar\zeta_\tau$ and $\bar\pi_\tau$:}
\nonumber
&\forall\widetilde\zeta\!\in\! L^\infty(\GC),\ \,
0\!\le\!\widetilde\zeta\!\le\!\bar\zeta_\tau(t):\ \
\\\label{e6:delam-ls-mix-semi-stab2}&
\qquad\qquad\quad
\calE\big(t,\bar u_\tau(t),\bar\zeta_\tau(t),\bar\pi_\tau(t)\big)\le\calE\big(t,u(t),\widetilde\zeta,\bar\pi_\tau(t)\big)
+\calR_0\big(\widetilde\zeta{-}\bar\zeta_\tau(t)\big),
\\
\nonumber&\forall\widetilde\pi\!\in\! H^1(\GC;\R^{d-1}):
\\\label{e6:delam-ls-mix-semi-stab1}
&\qquad\qquad\quad
\calE\big(t,\bar u_\tau(t),\underline\zeta_\tau(t),\bar\pi_\tau(t)\big)\le\calE\big(t,\bar u_\tau(t),\underline\zeta_\tau(t),\widetilde\pi\big)
+\calR_1\big(\widetilde\pi{-}\bar\pi_\tau(t)\big),
\intertext{and the energy (im)balance:}
\nonumber&
\calE\big(t_2,\bar u_\tau(t_2),\bar\zeta_\tau(t_2),\bar\pi_\tau(t_2)\big)
+\Diss_{\calR_1}^{}\big(\bar\pi_\tau;[t_1,t_2]\big)
\\[-.3em]\label{e6:delam-ls-mix-semi-engr}&\hspace{2em}
+\calR_0\big(\bar\zeta_\tau(t_2){-}\bar\zeta_\tau(t_1)\big)
\le
\calE\big(t_1,\bar u_\tau(t_1),\bar\zeta_\tau(t_1),\bar\pi_\tau(t_1)\big)
-\!\int_{t_1}^{t_2}\!\!
\big\langle\DT f_1,\underline u_\tau\big\rangle\,\d t,
\end{align}\end{subequations}
which is to hold for all $t\!\in\!I$
and for all $0\le t_1<t_2\le T$ of the form $t_i=k_i\tau$ for some $k_i\!\in\!\N$.
\end{proposition}

\noindent{\it Sketch of the proof.}
Writing optimality condition for \eqref{minimize-1} in terms of $u$,
one arrives at \eqref{e6:delam-mix-ls-VI}, and comparing the value of
\eqref{minimize-1} at $(u_\tau^k,\pi_\tau^k)$ with its value at
$(u_\tau^k,\widetilde\pi)$ and using the degree-1 homogeneity of $\calR_1$,
one arrives at \eqref{e6:delam-ls-mix-semi-stab1}.

Comparing the value of \eqref{minimize-2} at $\zeta_\tau^k$ with its
value at $\widetilde\zeta$ and using the degree-1 homogeneity of
$\calR_0$, one arrives at \eqref{e6:delam-ls-mix-semi-stab2}.

In obtaining \eqref{e6:delam-ls-mix-semi-engr},
we compare the value of \eqref{minimize-1} at the minimizer
$(u_\tau^k,\pi_\tau^k)$ with the value at $(u_\tau^{k-1},\pi_\tau^{k-1})$
and  the value of \eqref{minimize-2}  at the minimizer $\zeta_\tau^k$ with
the value at $\zeta_\tau^{k-1}$ and we benefit from the cancellation of the 
terms $\pm\calE(k\tau,u_\tau^k,\zeta_\tau^{k-1},\pi_\tau^k)$.
We also use the discrete by-part integration (=summation) for the $f_1$-term.

Then, using \eqref{e6:delam-ls-mix-semi-engr}
for $t_1=0$ and the coercivity of $\calE(t,\cdot,\cdot,\cdot)$ due to
the assumptions \eqref{ass}, we obtain also the a-priori estimates
\eqref{e6:delam-mix-ls-estimates}.
$\hfill\Box$

\medskip

The cancellation effect in the above proof is typical in fractional-step 
methods, cf.\ e.g.\ \cite[Remark~8.25]{Roub13NPDE} and for specific usage in
fracture mechanics also \cite{LaOrSu10ESRM}. Further, note that
\eqref{e6:delam-mix-ls} is of a similar form as \eqref{def-mix-LS} and is
thus prepared to make a limit passage for $\tau\to0$:

\begin{proposition}[Convergence towards local solutions]\label{prop-conv}
Let (\ref{ass}a-d) hold and $(\bar u_\tau,\bar\zeta_\tau,\bar\pi_\tau)$ be an 
approximate solution obtained by \eqref{delam-small-II-LS-min}.
Then, considering a sequence $\tau=\tau_n=T/n$ with $n\to\infty$, there exists 
a subsequence $\{(\bar u_\tau,\bar\zeta_\tau,\bar\pi_\tau)\}_{\tau>0}$
and $u\in{\rm B}(I;H^1(\Omega{\setminus}\GC;\R^d))$
with $\JUMP{u}{_{\rm N}}\ge0$ on $\SC$ and
$\zeta\in{\rm B}(I; L^\infty(\SC))\cap{\rm BV}(I;L^1(\SC))$ and
$\pi\in{\rm B}(I; H^1(\SC;\R^{d-1}))\cap{\rm BV}(I;L^1(\SC;\R^{d-1}))$
such that
\begin{subequations}\label{e4:delam-mix-semiimplicit-conv}
\begin{align}\label{e4:delam-mix-semiimplicit-conv-u}
&\bar u_\tau(t)\to u(t)&&\text{in } H^1(\Omega{\setminus}\GC;\R^d)
&&\hspace*{-3em}\text{for all }t\in I,&&
\\\label{e4:delam-mix-semiimplicit-conv-zeta}&
\bar\zeta_\tau(t)\weaks\zeta(t)&&\text{in } L^\infty(\GC),
&&\hspace*{-3em}\text{for all }t\in I,&&
\\\label{e4:delam-mix-semiimplicit-conv-pi}&
\bar\pi_\tau(t)\to\pi(t)&&\text{in } H^1(\GC;\R^{d-1})
&&\hspace*{-3em}\text{for all }t\in I.&&
\end{align}
\end{subequations}
Moreover, any $(u,\zeta,\pi)$ obtained by this way is a local solution
to the delamination problem in the sense of Definition~\ref{def-LS}.
\end{proposition}

\noindent{\it Proof.}
By Helly's selection principle \cite{Hell12LFO}, cf.\ also
e.g.\ \cite{Miel05ERIS,Miel11DEMF} for a more general version and usage
in rate-independent processes, we choose a subsequence and
$\zeta,\,\underline\zeta\in{\rm B}(I;L^\infty(\GC))\cap{\rm BV}(I;L^1(\GC))$
and $\pi\in{\rm B}(I;H^1(\GC;\R^{d-1}))\cap{\rm BV}(I;L^1(\GC;\R^{d-1}))$
so that
\begin{subequations}\label{eq3:LS-conv-to-z}\begin{align}
&\bar\zeta_\tau(t)\weak\zeta(t)&&\&&&\underline\zeta_\tau(t)\weak\underline\zeta(t)&&\text{ in }\ L^\infty(\GC)\ \text{ for all }\ t\!\in\!I,
\\
&\bar\pi_\tau(t)\weak\pi(t)&&
&&
&&\text{ in }\ H^1(\GC;\R^{d-1})\ \text{ for all }\ t\!\in\!I.
\end{align}\end{subequations}
Now, for a fixed $t\in I$, by Banach's selection principle, we select
(for a moment) further subsequence so that
\begin{align}\label{weak-conv-u}
\bar u_\tau(t)\weak u(t)\quad\text{ in }\ H^1(\Omega{\setminus}\GC;\R^d).
\end{align}
We further use that $\bar u_\tau(t)$ minimizes
$\calE(t_\tau,\cdot,\underline\zeta_\tau(t),\bar\pi_\tau)$
with $t_\tau:=\min\{k\tau\ge t;\ k\in\N\}$.
Obviously, $t_\tau\to t$ for $\tau\to0$
and, by the weak-lower-semicontinuity argument, we can easily
see that $u(t)$ minimizes the strictly convex functional
$\calE(t,\cdot,\underline\zeta(t),\pi(t))$.
Thus $u(t)$ is determined uniquely so that, in fact, we did not need to
make further selection of a subsequence, and this procedure
can be performed for any $t$  by using the same subsequence already
selected for \eqref{eq3:LS-conv-to-z}. Also,
$u:I\to H^1(\Omega{\setminus}\GC;\R^d)$ is measurable because
$\underline\zeta$ and
$\pi$ are measurable, and $\pl_u\calE(t,u(t),\underline\zeta(t),\pi(t))\!\ni\!0$
for all $t$.

The key ingredient is improvement of \eqref{weak-conv-u} for the strong convergence
of displacements: by using \eqref{e6:delam-mix-ls-VI} for $v=u(t)$ (which is a
legal test because the limit $u(t)$ satisfies the unilateral
constraint $\JUMP{u(t)}{_{\rm N}}\ge0$ on $\GC$), we have
\begin{align}\nonumber
&\hspace*{-.1em}\int_{\Omega{\setminus}\GC}\!\!\!
\mathbb{C}e(\bar u_\tau(t){-}u(t)){:}e(\bar u_\tau(t){-}u(t))\,\d x
\le\int_{\Omega{\setminus}\GC}\!\!\!\mathbb{C}e(\bar u_\tau(t){-}u(t)){:}
e(\bar u_\tau(t){-}u(t))\,\d x
\\\nonumber&\hspace*{5em}
+\int_{\GC}\!\!\!\underline\zeta_\tau(t)\Big(\kappa_{_{\rm N}}
\JUMP{\bar u_\tau(t){-}u(t)}{_{\rm N}}^{\!\!\!\!2}
+\kappa_{_{\rm T}}
\big|\JUMP{\bar u_\tau(t){-}u(t)}{_{\rm T}}\big|^2\Big)\d S
\\\nonumber&\qquad
\le\int_{\Omega{\setminus}\GC}\!\!\!
\mathbb{C} e(u(t)){:}e(u(t){-}\bar u_\tau(t))
\,\d x-\big\langle f_1(t_\tau),\bar u_\tau(t){-}u(t)\big\rangle
\\[-.3em]\nonumber
&\hspace*{5em}
+\int_{\GC}\underline\zeta_\tau(t)\Big(
\kappa_{_{\rm N}}\JUMP{u(t)}{_{\rm N}}{\cdot}\JUMP{u(t){-}\bar u_\tau(t)}{_{\rm N}}\!
\\[-.3em]\label{e6:delam-ls-strong}
&\hspace*{11em}
+\kappa_{_{\rm T}}\big(\JUMP{u(t)}{_{\rm T}}\!\!-\mathbb{T}\bar\pi_\tau(t)\big){\cdot}\JUMP{u(t){-}\bar u_\tau(t)}{_{\rm T}}\Big)\,
\d S\to0
\end{align}
with again $t_\tau:=\min\{k\tau\ge t;\ k\in\N\}$.
To prove this limit in \eqref{e6:delam-ls-strong} for $\tau\to0$,
we may simply use  $\JUMP{\bar u_\tau(t)}{}\weak\JUMP{u}{}$ in
$H^{1/2}(\GC;\R^d)$ so strongly in $L^2(\GC;\R^d)$ and
$\underline\zeta_\tau(t)\weaks\underline\zeta(t)$ in $L^\infty(\GC)$
so that
$\underline\zeta_\tau(t)\JUMP{\bar u_\tau(t)}{}\weak
\underline\zeta(t)\JUMP{u(t)}{}$
in $ L^2(\GC;\R^d)$. Due to the bound $\|\bar\pi_\tau(t)\|_{H^1(\GC;\R^{d-1})}$
and the compact embedding $H^1(\GC)\Subset L^2(\GC)$,
also $\bar\pi_\tau(t)\to\pi(t)$ in $L^2(\GC;\R^{d-1})$ and
thus $\underline\zeta_\tau(t)\mathbb{T}\bar\pi_\tau(t){}\weak
\underline\zeta(t)\mathbb{T}\pi(t)$ in $ L^2(\GC;\R^{d-1})$.
Then the convergence in \eqref{e6:delam-ls-strong} is trivial. We then obtain
 the strong convergence \eqref{e4:delam-mix-semiimplicit-conv-u}.

For the strong convergence \eqref{e4:delam-mix-semiimplicit-conv-pi}, we use
the information from the discrete flow-rule for $\pi$ obtained as an
optimality condition for \eqref{minimize-1} with respect to $\pi$, written as
\begin{align}\label{e6:delam-ls-mix-flow-for-pi-disc}
\bar{\mathfrak{f}}_\tau^{}\!
+\underline\zeta_\tau
\kappa_{_{\rm T}}\big(\JUMP{\bar u_\tau}{_{\rm T}}\!{-}\mathbb{T}\bar\pi_\tau\big)
+\kappa_{\scriptscriptstyle\rm{H}}\bar\pi_\tau
+\kappa_{_{\rm G}}{\rm div}_{\mbox{\tiny S}}\,\nablaS\bar\pi_\tau=0
\ \text{ with }\
\bar{\mathfrak{f}}_\tau\!\in\!N_{B_{\sigma_{\mathrm{yield}}}}(\DT\pi_\tau)
\end{align}
with $N_{B_{\sigma_{\mathrm{yield}}}}$ denoting the set-valued mapping
$\R^{d-1}\rightrightarrows\R^{d-1}$ defined as the normal cone
to the ball $B_{\sigma_{\mathrm{yield}}}\subset\R^{d-1}$ of the radius
$\sigma_{\mathrm{yield}}$
centered at the origin. The meaning of $\bar{\mathfrak{f}}_\tau$ is the discrete
driving force for the interfacial plasticity evolution. Fixing a time
instant $t$, we can thus assume $\bar{\mathfrak{f}}_\tau(t)$ bounded in
$ L^\infty(\GC;\R^{d-1})$ and use \eqref{e6:delam-ls-mix-flow-for-pi-disc} at
time $t$ tested by $\bar\pi_\tau(t)-\pi(t)$ to execute the limit passage
\begin{align}\nonumber
&\int_{\GC}\!\!\kappa_{_{\rm G}}|\nablaS\bar\pi_\tau(t)-\nablaS\pi(t)|^2\d S
\\[-.3em]\nonumber&\quad=\int_{\GC}\!\!\Big(\bar{\mathfrak{f}}_\tau(t)+\underline\zeta_\tau(t)
\kappa_{\scriptscriptstyle\rm{T}}\big|\JUMP{\bar u_\tau(t)}{_{\rm T}}\!{-}
\mathbb{T}\bar\pi_\tau(t)\big|^2\!
+\kappa_{\scriptscriptstyle\rm{H}}|\bar\pi_\tau(t)|^2\Big)
{\cdot}(\bar\pi_\tau(t){-}\pi(t))
\\[-.3em]\label{e6:delam-ls-strong+}&\hspace{16em}
-\kappa_{_{\rm G}}\nablaS\pi(t){:}\nablaS(\bar\pi_\tau(t){-}\pi(t))\d S\to0
\end{align}
where we again used the compact embedding $H^1(\GC)\Subset L^2(\GC)$.
Thus the strong convergence \eqref{e4:delam-mix-semiimplicit-conv-pi} follows.

The BV-functions (here in particular both BV-functions $\zeta(\cdot)$ and
$\underline\zeta(\cdot)$) are continuous
everywhere except at most countable number of times, let us denote
this set of jumps by $J$. Then we have $\zeta(t)=\underline\zeta(t)$
for any $t\in I{\setminus}J$.
In particular, $\partial_u\calE(t,u(t),\zeta(t),\pi(t))=
\partial_u\calE(t,u(t),\underline\zeta(t),\pi(t))\ni0$
for such $t$, which proves \eqref{e6:delam-mix-ls-VI}.

Now we can already pass to the limit in \eqref{e6:delam-mix-ls}. The limit
passage in \eqref{e6:delam-mix-ls-VI} for all $t\in I{\setminus}J$
simple just by continuity; note that we need $\zeta(t)=\underline\zeta(t)$
for all $t$ except from $J$. Thus \eqref{def-mix-ls-VI} is obtained.

For  the limit passage in  the semi-stability
\eqref{e6:delam-ls-mix-semi-stab2}
 towards \eqref{def-ls-mix-semi-stab2}, we use the
so-called mutual recovery sequence
\begin{align}\label{eq5:recov-seq-semistability}
   \widetilde\zeta_\tau(x):=\begin{cases}\
\bar\zeta_\tau(t,x)\widetilde\zeta(x)/\zeta(t,x) &
\text{if }\zeta(t,x)>0,
\\\ 0 & \text{if }\zeta(t,x)=0.\end{cases}
\end{align}
with $0\le\widetilde\zeta\le\zeta(t)$ given.
After substituting $\widetilde\zeta_\tau$ in place of $\widetilde\zeta$
into \eqref{e6:delam-ls-mix-semi-stab2}, we can easily pass to
\eqref{def-ls-mix-semi-stab2} by continuity, namely
\begin{align}\nonumber
0&\le\lim_{\tau\to0}\int_{\GC}\!\!\big(\!\!\!\!\!\!\!\!\!
\ddd{\widetilde\zeta_\tau{-}\bar\zeta_\tau(t)}
{converges}{weakly* in $L^\infty(\GC)$}\!\!\!\!\!\!\big)
\Big(\!
\ddd{\frac{\kappa_{_{\rm N}}}2\JUMP{\bar u_\tau(t)}{_{\rm N}}^{\!\!\!\!2}+
\frac{\kappa_{_{\rm T}}}2\big|\JUMP{\bar u_\tau(t)}{_{\rm T}}\!{-}\mathbb{T}\bar\pi_\tau(t)\big|^2}{converges in $H^{1/2}(\GC)\subset L^1(\GC)$}{}\!\!\!
-a_{_{\rm I}}\Big)\d S\!
\\\label{limit-semistab-zeta}&=
\int_{\GC}\!\!(\widetilde\zeta{-}\zeta(t))
\Big(\frac{\kappa_{_{\rm N}}}2\JUMP{u(t)}{_{\rm N}}^{\!\!\!\!2}+
\frac{\kappa_{_{\rm T}}}2\big|\JUMP{u(t)}{_{\rm T}}\!{-}\mathbb{T}\pi(t)\big|^2
-a_{_{\rm I}}\Big)\d S\!
\end{align}
which is just the semistability \eqref{def-ls-mix-semi-stab2}.

It is important
that $0\le\widetilde\zeta_\tau\le\bar\zeta_\tau(t)$
a.e.\ on $\GC$ and, since $\bar\zeta_\tau(t)\weaks\zeta(t)$,
also $\widetilde\zeta_\tau\weaks\widetilde\zeta$ in $L^\infty(\GC)$.
 For this explicit construction \eqref{eq5:recov-seq-semistability},
cf.~also \cite[Lemma 6.1]{MaiMie05EREM} or \cite[Formula (3.71)]{RoScZa09QDP}.

Also the limit passage in \eqref{e6:delam-ls-mix-semi-stab1} towards
\eqref{def-ls-mix-semi-stab1} is simple just by continuity 
because we already proved the strong convergence
\eqref{e4:delam-mix-semiimplicit-conv-pi} otherwise the weak convergence
would serve here too by semi-continuity arguments. The
mutual recovery sequence can be even taken simply constant, namely
$\widetilde\pi_\tau=\widetilde\pi$, so that:
\begin{align}\nonumber
&\int_{\GC}\!\!\underline\zeta(t)\frac{\kappa_{_{\rm T}}}2
\big|\JUMP{u(t)}{_{\rm T}}\!{-}\mathbb{T}\pi(t)\big|^2\!
+\frac{\kappa_{\scriptscriptstyle\rm H}}2|\pi(t)|^2\!
+\frac{\kappa_{_{\rm G}}}2|\nablaS\pi(t)|^2\d S
\\\nonumber
&\ \ =\lim_{\tau\to0}
\int_{\GC}\!\!\underline\zeta_\tau(t)\frac{\kappa_{_{\rm T}}}2
\big|\JUMP{\bar u_\tau(t)}{_{\rm T}}\!{-}\mathbb{T}\bar\pi_\tau(t)\big|^2\!
+\frac{\kappa_{\scriptscriptstyle\rm H}}2|\bar\pi_\tau(t)|^2\!
+\frac{\kappa_{_{\rm G}}}2|\nablaS\bar\pi_\tau(t)|^2\d S
\\\nonumber
&\ \ \le\lim_{\tau\to0}\int_{\GC}\!\underline\zeta_\tau(t)
\frac{\kappa_{_{\rm T}}}2\big|\JUMP{\bar u_\tau(t)}{_{\rm T}}\!{-}\mathbb{T}\widetilde\pi\big|^2
+\frac{\kappa_{\scriptscriptstyle\rm H}}2|\widetilde\pi|^2
+\frac{\kappa_{_{\rm G}}}2|\nablaS\widetilde\pi|^2\!
+\sigma_{\mathrm{yield}}\big|\widetilde\pi{-}\bar\pi_\tau(t)\big|\d S
\\\label{e6:delam-ls-mix-semi-stab1-disc}
&\ \ =
\int_{\GC}\!\underline\zeta(t)
\frac{\kappa_{_{\rm T}}}2\big|\JUMP{u(t)}{_{\rm T}}\!{-}\mathbb{T}\widetilde\pi\big|^2
+\frac{\kappa_{\scriptscriptstyle\rm H}}2|\widetilde\pi|^2
+\frac{\kappa_{_{\rm G}}}2|\nablaS\widetilde\pi|^2\!
+\sigma_{\mathrm{yield}}\big|\widetilde\pi{-}\pi(t)\big|\d S.
\end{align}

The limit passage in  the energy (im)balance 
\eqref{e6:delam-ls-mix-semi-engr} towards \eqref{def-ls-mix-engr} relies on the
(strong$\times$weak$\times$strong)-continuity of
$\calE(t,\cdot,\cdot,\cdot)$ on its definition domain.
First we need to extend \eqref{e6:delam-ls-mix-semi-engr}
for all $t_1$ and $t_2$. By \eqref{E-delam-small-k-II},
we have $\pl_t\calE\big(t,\underline u_\tau,\underline\zeta_\tau,\underline \pi_\tau\big)
=\langle\DT f_1,\underline u_\tau\rangle$, and by the assumption
$f_1\in W^{1,1}(I; H^1(\Omega{\setminus}\GC))$, it is easy to see that
\begin{align}\nonumber
&\calE\big(t_2,\bar u_\tau(t_2),\bar\zeta_\tau(t_2),\bar\pi_\tau(t_2)\big)
+\calR_0\big(\bar\zeta_\tau(t_2){-}\bar\zeta_\tau(t_1)\big)
+\Diss_{\calR_1}^{}\big(\bar\pi_\tau;[t_1,t_2]\big)
\\[-.2em]&\label{e6:delam-ls-engr+}
\hspace{6em}
\le\calE\big(t_1,\bar u_\tau(t_1),\bar\zeta_\tau(t_1),\bar\pi_\tau(t_1)\big)
-\!\int_{t_1}^{t_2}\!\!
\big\langle\DT f_1,\underline u_\tau\big\rangle
\,\d t+\mathscr{O}(\tau),
\end{align}
for all $0\le t_1\le t_2\le T$. By \eqref{e4:delam-mix-semiimplicit-conv}
and by the arguments we already used for \eqref{e6:delam-ls-strong},
we can easily see that $\calE(t,\bar u_\tau(t),\bar\zeta_\tau(t),\bar\pi_\tau(t))
\to\calE(t,u(t),\zeta(t),\pi(t))$, which is to be used for
\eqref{e6:delam-ls-engr+} both for $t=t_1$ and $t=t_2$.
$\hfill\Box$

\begin{remark}[{A}pproximate  maximum-dissipation principle]\label{rem-DMP}
\upshape
One can devise the discrete analog of the integrated maximum-dissipation
principle \eqref{IMDP} straightforwardly for the left-continuous
interpolants \eqref{e3:def-of-interpolants},
required however to hold only asymptotically.
More specifically, in analog to \eqref{IMDP} formulated
equivalently for all $[0,t]$ instead of $[t_1,t_2]$, one can expect
an {\it Approximate Maximum-Dissipation Principle} (AMDP) in the form
\begin{subequations}\label{AMDP}
\begin{align}\label{AMDP-a}
&\!\int_0^t\!\bar{\mathfrak f}_\tau\,\d\bar\pi_\tau
\ \stackrel{\mbox{\bf ?}}{\sim}\ \Diss_{\calR_1}(\bar\pi_\tau;[0,t])
\;\ \ \ \text{ for some }\ \ \bar{\mathfrak f}_\tau
\!\in\!-\partial_\pi\bar{\calE}_\tau(\cdot,\bar u_\tau,\underline\zeta_\tau,\bar\pi_\tau),
\\\label{AMDP-b}
&\!\int_0^t\!\bar{\mathfrak g}_\tau\,\d\bar\zeta_\tau
\ \stackrel{\mbox{\bf ?}}{\sim}\ \Diss_{\calR_0}(\bar\zeta_\tau;[0,t])
\ \ \ \text{ for some }\ \ \bar{\mathfrak g}_\tau
\!\in\!
-\partial_\zeta\bar{\calE}_\tau(\cdot,\bar u_\tau,\bar\zeta_\tau,\bar\pi_\tau),
\end{align}\end{subequations}
where again the integrals are the lower Riemann-Stieltjes integrals as in
\eqref{IMDP} and where $\bar{\calE}_\tau(\cdot,u,\zeta,\pi)$ is the
left-continuous piecewise-constant interpolant of the values
$\calE(k\tau,u,\zeta,\pi)$, $k=0,1,...,T/\tau$.
Moreover, ''${\stackrel{\mbox{\footnotesize\bf ?}}{\sim}}$'' in \eqref{AMDP}
means that the equality holds possibly only asymptotically for $\tau\to0$ but
even this is rather only desirable and not always valid.
Anyhow, loadings which, under given geometry of the specimen,
lead to  rate-independent  slides  where the solution
is absolutely continuous will always comply with AMDP \eqref{AMDP}.
Also, some finite-dimensional examples  of ``damageable springs''
in \cite{MieRou15RIST,Roub??MDLS} show that this AMDP can detect too
early rupturing local solutions (in particular the energetic ones)
while it generically holds for solutions obtained by the
algorithm \eqref{delam-small-II-LS-min}.
In our model, too early rupturing may also mean unphysical sliding
into less dissipative Mode I even in situations when clearly Mode II
should be active, cf.\ also the computational experiments in
\cite{VoMaRo??BEIM}. Generally speaking, \eqref{AMDP}
should rather be a-posteriori checked to justify the
(otherwise not physically based) simple and numerically efficient
fractional-step-type semi-implicit algorithm \eqref{delam-small-II-LS-min}
from the perspective of the stress-driven solutions in particular
situations and possibly to provide a valuable information
that can be exploited to adapt time or space discretisation towards
better accuracy in \eqref{AMDP} and thus close towards the stress-driven
scenario. Actually, for the piecewise-constant interpolants, we can simply 
evaluate the integrals explicitly, so that AMDP \eqref{AMDP} reads
\begin{subequations}\label{e1:AMDP}\begin{align}
&
\sum_{k=1}^K\int_{\GC}\!\!\mathfrak f_\tau^{k-1}(\pi_\tau^k-\pi_\tau^{k-1})\,\d S
\ \stackrel{\mbox{\large\bf ?}}{\sim}\
\sum_{k=1}^K\int_{\GC}\!\!
\sigma_\mathrm{yield}^{}\big|\pi_\tau^k{-}\pi_\tau^{k-1}\big|\,\d S
\ \ \ \ \text{ and}
\\[-.3em]&\!
\sum_{k=1}^K\int_{\GC}\!\!\mathfrak g_\tau^{k-1}(\zeta_\tau^k{-}\zeta_\tau^{k-1})\,\d S
\ \stackrel{\mbox{\large\bf ?}}{\sim}\
\int_{\GC}\!\!a_{_{\rm I}}\big(\zeta_0{-}\zeta_\tau^K\big)\,\d S
\\&\nonumber\qquad\qquad\text{ where }
\ K=\max\{k\!\in\!\mathbb N;\ k\tau\le t\}\ \ \text{ and }
\\&\nonumber\qquad\qquad\text{ where }\
\mathfrak f_\tau^k\in-\partial_\pi\calE(u_\tau^k,\zeta_\tau^{k-1},\pi_\tau^k)
\ \ \,\text{ and}\ \ \,
\mathfrak g_\tau^k\in-\partial_\zeta\calE(u_\tau^k,\zeta_\tau^k,\pi_\tau^k).
\end{align}\end{subequations}
Always, the left-hand sides in \eqref{e1:AMDP} are below the right-hand
sides, and one can a-posteriori check the residua depending on $t$ (or
possibly also on space, cf.\ \cite{VoMaRo??BEIM}).
\end{remark}

\section{Numerical approximation and computational experiments}\label{sect-BEM}
%        ~~~~~~~~~~~~~~~~~~~~~~~~~~~~~~~~~~~~~~~~~~~~~~~~~~~~~

Let us assume $\Omega\subset\R^d$ to be a polyhedral domain with
$\GD,\GN,\GC\subset\R^{d-1}$ also polyhedral. We
outline briefly the discretisation by the finite-element method.
In the simplest variant, $\Omega$ is discretised
by a triangular mesh $\mathscr{T}_h$ consistently with the
boundaries $\GD$ and $\GC$ with $h>0$ denoting the mesh parameter,
and the polynomial P1-elements for $u$, P0-elements for $\zeta$, and
$P1$-elements for $\pi$ are employed.
Applying such an approximation to \eqref{delam-small-II-LS-min},
we thus arrive at two linear-quadratic programming problems:
\begin{subequations}\label{LQ}
\begin{align}\label{LQ-1}
&\left.\begin{array}{ll}
\text{minimize}&\calE(k\tau,u,\zeta_{\tau h}^{k-1},\pi)
+\calR_1(\pi{-}\pi_{\tau h}^{k-1})
\\[.2em]\text{subject to}&(u,\pi)\in H^1(\Omega{\setminus}\GC;\R^d)\times H^1(\GC;\R^{d-1}),
\\[.2em] &(u,\pi)\text{ element-wise linear on }\mathscr{T}_h,
\end{array}\right\}
\intertext{and, denoting the unique solution as $(u_{\tau h}^k,\pi_{\tau h}^k)$,}
\label{LQ-2}
&\left.\begin{array}{ll}
\text{minimize}&\calE(k\tau,u_{\tau h}^k,\zeta,\pi_{\tau h}^k)
+\calR_0(\zeta{-}\zeta_{\tau h}^{k-1})
\\[.2em]\text{subject to}&\zeta\in L^\infty(\GC),\ \ 0\le\zeta\le1,
\\[.2em]&\zeta\text{ element-wise constant on }\mathscr{T}_h,
\end{array}\ \ \ \ \ \ \ \,\right\}
\end{align}
\end{subequations}
and denote its (possibly not unique) solution by $\zeta_{\tau h}^k$.
Existence of such finite-dimensional solutions
$(u_{\tau h}^k,\zeta_{\tau h}^k,\pi_{\tau h}^k)$
is even simpler than in Section~\ref{sec-disc} because the considered
linear spaces are finite-dimensional. Numerically, the solution can be obtained
non-iteratively after a finite-number of steps if the linear-quadratic solver
used for \eqref{LQ} is implemented in this way.
 More in detail, $\calR_1$ in the cost functional in \eqref{LQ-1}
is nonsmooth and, only after applying the Mosco-type transformation
as e.g.\ in \cite[Lemma~4]{Roub02EMMP}, one obtains
truly a quadratic programming problem (QP) if $d=2$ or
a so-called  second-order cone programming  problem (SOCP)
if $d=3$ when $\calR_1$ does not have a polyhedral graph, cf.\ e.g.\
\cite{AliGol03SOCP,Stur02IIPM} for the SOCP algorithms.

\begin{proposition}[{U}nconditional convergence towards local solutions]
\label{cor-convergence}
Let again (\ref{ass}a-d) holds  and let the spatial discretisation
refines everywhere, i.e.\
$\lim_{h\to0}^{}\sup_{\triangle\in\mathscr{T}_h}\mathrm{diam}(\triangle)=0$. 
Then, the solution to the recursive alternating-minimization problem \eqref{LQ}
exists and is numerically stable, i.e., in terms of the time-interpolants,
\begin{subequations}\label{estimates}
\begin{align}\label{estimate1}
&\big\|\bar u_{\tau h}\big\|_{ L^\infty(I; H^1(\Omega{\setminus}\GC;\R^d))}
\le C,&&
\\\label{estimate2}
&\big\|\bar\zeta_{\tau h}\big\|_{ L^\infty(\SC)
\cap{\rm BV}(I; L^1(\GC))}\le C,
\\\label{estimate3}
&\big\|\bar\pi_{\tau h}\big\|_{ L^\infty(I; H^1(\GC;\R^{d-1}))
\cap{\rm BV}(I; L^1(\GC;\R^{d-1}))}\le C
&&
\end{align}\end{subequations}
with some $C$ independent of $\tau>0$ and $h>0$.
This solution satisfies the analog of \eqref{e6:delam-mix-ls}
with the test functions $\widetilde u$, $\widetilde\zeta$, and $\widetilde\pi$
ranging over the above specified FEM-subspaces.
Moreover, if $\tau\to0$ and
$h\to0$, then in terms of subsequences, like in Proposition~\ref{prop-conv},
it converges to local solutions to the delamination problem
\eqref{Gm+}--\eqref{IC}:
\begin{subequations}\label{e4:delam-mix-conv}
\begin{align}\label{e4:delam-mix-conv-u}
&\bar u_{\tau h}(t)\to u(t)&&\text{in }H^1(\Omega{\setminus}\GC;\R^d)
&&\hspace*{-3em}\text{for all }t\in I,&&
\\\label{e4:delam-mix-conv-zeta}&
\bar\zeta_{\tau h}(t)\weaks\zeta(t)&&\text{in }L^\infty(\GC)
&&\hspace*{-3em}\text{for all }t\in I,&&
\\\label{e4:delam-mix-conv-pi}&
\bar\pi_{\tau h}(t)\to\pi(t)&&\text{in }H^1(\GC;\R^{d-1})
&&\hspace*{-3em}\text{for all }t\in I.&&
\end{align}
\end{subequations}
\end{proposition}

\noindent{\it Sketch of the proof.}
The arguments of the proof of Proposition~\ref{prop-conv} can be
applied with only slight and mostly straightforward variation. 
Let us only briefly sketch differences 
beside that, of course, everywhere ``$\tau h$''
is written in place of the subscript ``$\tau$'' except in $t_\tau$.

The selection of converging subsequences is like in
\eqref{eq3:LS-conv-to-z}--\eqref{weak-conv-u}. Then, 
in \eqref{e6:delam-ls-strong}, one must use an element-wise affine approximant
of $u$ rather than directly $u$ itself. More in detail,
as $\bar u_{\tau h}(t)-u(t)$ is not a legal test function for the
Galerkin analog of \eqref{e6:delam-mix-ls-VI}, the estimate
\eqref{e6:delam-ls-strong} written with ``$\tau h$'' in place of
``$\tau$'' now modifies as
\begin{align}\nonumber
&\hspace*{-.1em}\int_{\Omega{\setminus}\GC}\!\!\!\!\!
\mathbb{C}e(\bar u_{\tau h}(t){-}u(t)){:}e(\bar u_{\tau h}(t){-}u(t))\,\d x
\\\nonumber&\quad
=\int_{\Omega{\setminus}\GC}\!\!\!\!\!\!\mathbb{C}e(\bar u_{\tau h}(t){-}u(t)){:}
e(\bar u_{\tau h}(t){-}\widetilde u_{t,h})
+\mathbb{C}e(\bar u_{\tau h}(t){-}u(t)){:}e(\widetilde u_{t,h}){-}u(t))\,\d x
\\\nonumber&\quad
\le\int_{\Omega{\setminus}\GC}\!\!\!\!\!\!\mathbb{C}e(\bar u_{\tau h}(t){-}u(t)){:}
e(\bar u_{\tau h}(t){-}\widetilde u_{t,h})
+\mathbb{C}e(\bar u_{\tau h}(t){-}u(t)){:}e(\widetilde u_{t,h}){-}u(t))\,\d x
\\[-.4em]\nonumber&\hspace*{5em}
+\int_{\GC}\!\!\!\underline\zeta_{\tau h}(t)\Big(\kappa_{_{\rm N}}
\JUMP{\bar u_{\tau h}(t){-}\widetilde u_{t,h}}{_{\rm N}}^{\!\!\!\!2}
+\kappa_{_{\rm T}}
\big|\JUMP{\bar u_{\tau h}(t){-}\widetilde u_{t,h}}{_{\rm T}}\big|^2\Big)\d S
\\\nonumber&\quad
\le\int_{\Omega{\setminus}\GC}\!\!\!
\mathbb{C} e(u(t)){:}e(\widetilde u_{t,h}{-}\bar u_{\tau h}(t))
+\mathbb{C}e(\bar u_{\tau h}(t){-}u(t)){:}e(\widetilde u_{t,h}){-}u(t))
\,\d x
\\[-.3em]\nonumber&\hspace*{2em}
-\big\langle f_1(t_\tau),\bar u_{\tau h}(t){-}\widetilde u_{t,h}\big\rangle
+\int_{\GC}\underline\zeta_{\tau h}(t)\Big(
\kappa_{_{\rm N}}\JUMP{\widetilde u_{t,h}}{_{\rm N}}{\cdot}\JUMP{\widetilde u_{t,h}{-}\bar u_{\tau h}(t)}{_{\rm N}}\!
\\[-.3em]\label{e6:delam-ls-strong-disc}
&\hspace*{10em}
+\kappa_{_{\rm T}}\big(\JUMP{\widetilde u_{t,h}}{_{\rm T}}\!\!-\mathbb{T}\bar\pi_{\tau h}(t)\big){\cdot}\JUMP{\widetilde u_{t,h}{-}\bar u_{\tau h}(t)}{_{\rm T}}\Big)\,
\d S\to0
\end{align}
where $\widetilde u_{t,h}$ is element-wise linear on $\mathscr{T}_h$ and
approximates $u(t)$ in the sense $\widetilde u_{t,h}\to u(t)$ in
$H^1(\Omega{\setminus}\GC;\R^d)$; such $\widetilde u_{t,h}$ always exists
provided only $h\to0$ because the spatial discretisation is supposed to
refine everywhere, and the possible dependence on
the rate of approximation of $u(t)$ on $t$ is unimportant for
\eqref{e6:delam-ls-strong-disc}.

Similarly, in \eqref{e6:delam-ls-strong+}, one must use an element-wise affine
approximant of $\pi$ rather than directly $\pi$ itself.
 More in detail,
\eqref{e6:delam-ls-strong+} written with ``$\tau h$'' in place of
``$\tau$'' modifies as
\begin{align}\nonumber
&\!\!\!\!
\int_{\GC}\!\!\kappa_{_{\rm G}}|\nablaS\bar\pi_{\tau h}(t)-\nablaS\pi(t)|^2\d S
=\int_{\GC}\!\!\kappa_{_{\rm G}}\!\nablaS(\bar\pi_{\tau h}(t){-}\pi(t)){:}
\nablaS(\bar\pi_{\tau h}(t){-}\widetilde\pi_{t,h})
\\[-.6em]\nonumber&\hspace{15em}
+\kappa_{_{\rm G}}\!\nablaS(\bar\pi_{\tau h}(t){-}\pi(t)){:}
\nablaS(\widetilde\pi_{t,h}{-}\pi(t))\,\d S
\\[-.1em]\nonumber&
=\!\!\int_{\GC}\!\!\!\Big(\bar{\mathfrak{f}}_{\tau h}(t)+\underline\zeta_{\tau h}(t)
\kappa_{\scriptscriptstyle\rm{T}}\big|\JUMP{\bar u_{\tau h}(t)}{_{\rm T}}\!{-}
\mathbb{T}\bar\pi_{\tau h}(t)\big|^2\!
+\kappa_{\scriptscriptstyle\rm{H}}|\bar\pi_{\tau h}(t)|^2\Big)
{\cdot}(\bar\pi_{\tau h}(t){-}\widetilde\pi_{t,h})
\\[-.3em]\label{e6:delam-ls-strong+disc}&\hspace{1em}
-\kappa_{_{\rm G}}\!\nablaS\pi(t){:}\nablaS(\bar\pi_{\tau h}(t){-}\widetilde\pi_{t,h})+\kappa_{_{\rm G}}\!\nablaS(\bar\pi_{\tau h}(t){-}\pi(t)){:}
\nablaS(\widetilde\pi_{t,h}{-}\pi(t))\d S\to0
\end{align}
where $\bar{\mathfrak{f}}_{\tau h}(t)$ is the discrete driving force
analogous as in \eqref{e6:delam-ls-mix-flow-for-pi-disc} and
again bounded in $L^\infty(\GC;\R^{d-1})$, and where
$\widetilde\pi_{t,h}$  is element-wise affine on $\mathscr{T}_h$ and
approximates $\pi(t)$ in the sense $\widetilde\pi_{t,h}\to\pi(t)$ in
$H^1(\GC;\R^{d-1})$; such $\widetilde\pi_{t,h}$ always exists
provided only $h\to0$ as the spatial discretisation
refines everywhere, and again the possible dependence on
the rate of approximation of $\pi(t)$ on $t$ is unimportant for
\eqref{e6:delam-ls-strong+disc}.

Instead of \eqref{eq5:recov-seq-semistability}, one can use the mutual 
recovery sequence:
\begin{align}\label{eq5:recov-seq-semistab}
   \widetilde\zeta_{\tau h}(x):=\begin{cases}\
\bar\zeta_{\tau h}(t,x)\big[\Pi_h^{(0)}\big(\widetilde\zeta/\zeta(t)\big)\big](x) &
\text{if }\big[\Pi_h\zeta(t)\big](x)>0,
\\\ 0 & \text{if }\big[\Pi_h\zeta(t)\big](x)=0\end{cases}
\end{align}
with $\Pi_h^{(0)}$ denoting the element-wise constant interpolation on $\GC$,
cf.~also \cite[Formula (4.35)]{MieRou09NARI}.
 If $z(x)=0$, then also $\widetilde z(x)=0$ because always
$0\le\widetilde z\le z$ and the fraction in \eqref{eq5:recov-seq-semistab}
can be defined arbitrarily  and valued in  $[0,1]$. The product of element-wise
constant functions $z_h$ and $\Pi_h^{(0)}(\widetilde z/z)$ is again
element-wise constant, hence $z_h\in Z_h$. As
$0\le\Pi_h^{(0)}(\widetilde z/z)\le1$, we have also
$0\le\widetilde z_h\le z_h$, hence $\widetilde z_h\in Z_h$ and
$\calR_0(\widetilde z_h{-}z_h)<\infty$. As
$\Pi_h^{(0)}(\widetilde z/z)\to\widetilde z/z$
strongly in any $L^p_{}(\GC)$, $p<+\infty$, and $z_h\weaks z$; here again
we rely on that the spatial discretisation is supposed to
refine everywhere. {F}rom \eqref{eq5:recov-seq-semistab} we have
$\widetilde z_h\weaks z(\widetilde z/z)=\widetilde z$ in fact in
$L^\infty_{}(\GC)$ due to the a priori bound of values in [0,1].
The limit passage from the discretised analog of
\eqref{e6:delam-ls-mix-semi-stab2} towards the semistability
\eqref{def-ls-mix-semi-stab2} is then completely
analogous to \eqref{limit-semistab-zeta}.

Also for the limit passage
in the spatially-discretised analog of \eqref{e6:delam-ls-mix-semi-stab1},
instead of just $\widetilde\pi$ fixed, one must use
\begin{align}
\widetilde\pi_{\tau h}:=\pi_{\tau h}-\Pi_h^{(1)}(\widetilde\pi{-}\pi)
\end{align}
with $\Pi_h^{(1)}$ denoting the element-wise affine interpolation on $\GC$,
cf.~also \cite[Formula (3.31)]{MieRou09NARI}.  A
modification of \eqref{e6:delam-ls-mix-semi-stab1-disc} is then
straightforward because $\widetilde\pi_{\tau h}\to\widetilde\pi$
strongly in $H^1(\GC;\R^{d-1})$; also here  
we rely on that the spatial discretisation is supposed to
refine everywhere. 
$\hfill\Box$

\medskip

The approximate  maximum-dissipation principle \eqref{e1:AMDP}
now reads as: 

\begin{subequations}\label{e1:AMDP+}\begin{align}\label{e1:AMDP+a}
&
\sum_{k=1}^K\int_{\GC}\!\!\mathfrak f_{\tau h}^{k-1}(\pi_{\tau h}^k-\pi_{\tau h}^{k-1})\,\d S
\ \ \stackrel{\mbox{\large\bf ?}}{\sim}\ \
\sum_{k=1}^K\int_{\GC}\!\!
\sigma_\mathrm{yield}^{}\big|\pi_{\tau h}^k{-}\pi_{\tau h}^{k-1}\big|\,\d S
\ \ \ \ \text{ and}
\\[-.3em]&\label{e1:AMDP+b}\!
\sum_{k=1}^K\int_{\GC}\!\!\mathfrak g_{\tau h}^{k-1}(\zeta_{\tau h}^k{-}\zeta_{\tau h}^{k-1})\,\d S
\ \ \stackrel{\mbox{\large\bf ?}}{\sim}\ \ 
\int_{\GC}\!\!a_{_{\rm I}}\big(\zeta_{\tau h}^K{-}\zeta_0\big)\,\d S
\end{align}\end{subequations}
where $\mathfrak f_{\tau h}^k\!\in\!-\partial_\pi\calE(u_{\tau h}^k,\zeta_{\tau h}^{k-1},\pi_{\tau h}^k)$ and $\mathfrak g_{\tau h}^k\!\in\!-\partial_\zeta\calE(u_{\tau h}^k,\zeta_{\tau h}^k,\pi_{\tau h}^k)$, and $K$ is as in \eqref{e1:AMDP}.

It is a noteworthy attribute of our problem that all inelastic processes occur
on the boundary $\GC$ while in the bulk domains $\Omega_1$ and $\Omega_2$
it is linear. This allows for elimination of nodal values inside $\Omega_1$
and $\Omega_2$ and considerable reduction of degrees of freedom by
considering only nodal or element values on $\GC$.

In fact, this idea has been systematically exploited even on the
continuous level when implementing the boundary-element method, cf.\
\cite{PaMaRo13BEMI,RoMaPa13QMMD,RoPaMA13QACV,TMGP11BEMA,VoMaRo??BEIM},
although it is still not fully supported by a convergence analysis like
Corollary~\ref{cor-convergence} due to general substantial
theoretical difficulties related to this method.

Anyhow, for the computational experiments presented here with the
goal to document rather modelling issues, we use a shortcut in
implementing the spatial discretisation \eqref{LQ} by exploiting
the collocation boundary-element method. Another numerical shortcut
was neglecting the gradient term by putting $\kappa_{_{\rm G}}=0$.

\begin{my-picture}{.95}{.2}{fig_m1}
\psfrag{GN}{\footnotesize $\GN$}
\psfrag{GD}{\footnotesize $\GD$}
\psfrag{GC}{\footnotesize $\GC$}
\psfrag{elastic}{\footnotesize elastic body}
\psfrag{obstacle}{\footnotesize rigid obstacle}
\psfrag{adhesive}{\footnotesize adhesive}
\psfrag{LC}{$L_c$}
\psfrag{L}{\footnotesize $L=\ $250\,mm}
\psfrag{H}{\footnotesize $H=$}
\psfrag{12.5}{\scriptsize 12.5\,mm}
\psfrag{loading}{\footnotesize loading}
\hspace*{-.5em}\vspace*{-.1em}{\includegraphics[width=.95\textwidth]{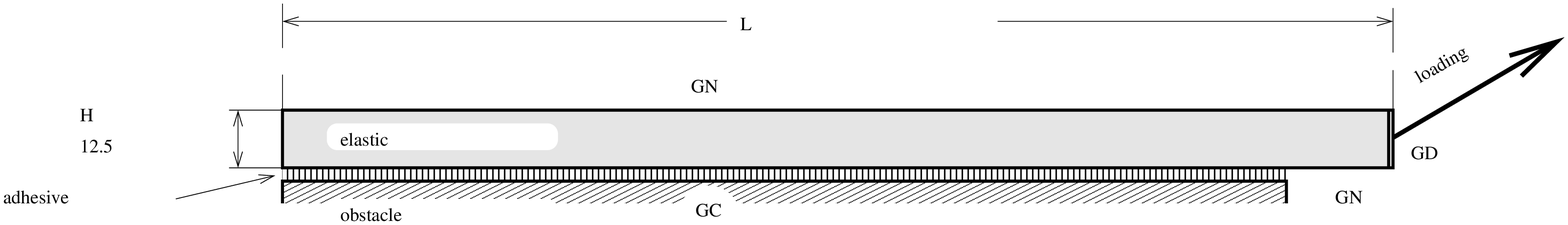}}
\end{my-picture}
\vspace*{-1.em}
\begin{center}
{\small Fig.\,\ref{fig_m1}.\ }
\begin{minipage}[t]{.70\textwidth}\baselineskip=8pt
{\small
Geometry and boundary conditions of the problem considered.
The length of the initially glued part $\GC$ is $0.9L=225\,$mm,
the adhesive layer has zero thickness.}
\end{minipage}
\end{center}

We demonstrate varying mode-mixity of delamination on
a relatively simple example motivated by the pull-push shear
experimental test used in engineering practice \cite{CoCa11ES}.
Intentionally, we use the same geometry, shown in Fig.~\ref{fig_m1}, as in \cite{RoMaPa13QMMD} in order to compare our maximally-dissipative
local solution with the energetic solution presented in \cite{RoMaPa13QMMD}.
In contrast to Sections~\ref{sect_AssocMod}--\ref{sec-disc},
only one bulk domain is considered and $\GC$ is a part of its
boundary but this modification is straightforward;
alternatively, one may also think about
$\Omega_2$ as a completely rigid body in the previous setting.
Here $\Omega_1$ is a two-dimensional rectangular domain glued on the
most of its bottom side $\GC$ with
the Dirichlet loading acting on the right-hand side $\GD$ in
the direction $(1,0.6)$, cf.\  Fig.~\ref{fig_m1}, increasing linearly in time
with velocity $1\,$mm/s.

\begin{my-picture}{.95}{.5}{fig_m2}
\hspace*{1.5em}\vspace*{-.1em}{\includegraphics[width=.8\textwidth,height=.45\textwidth]{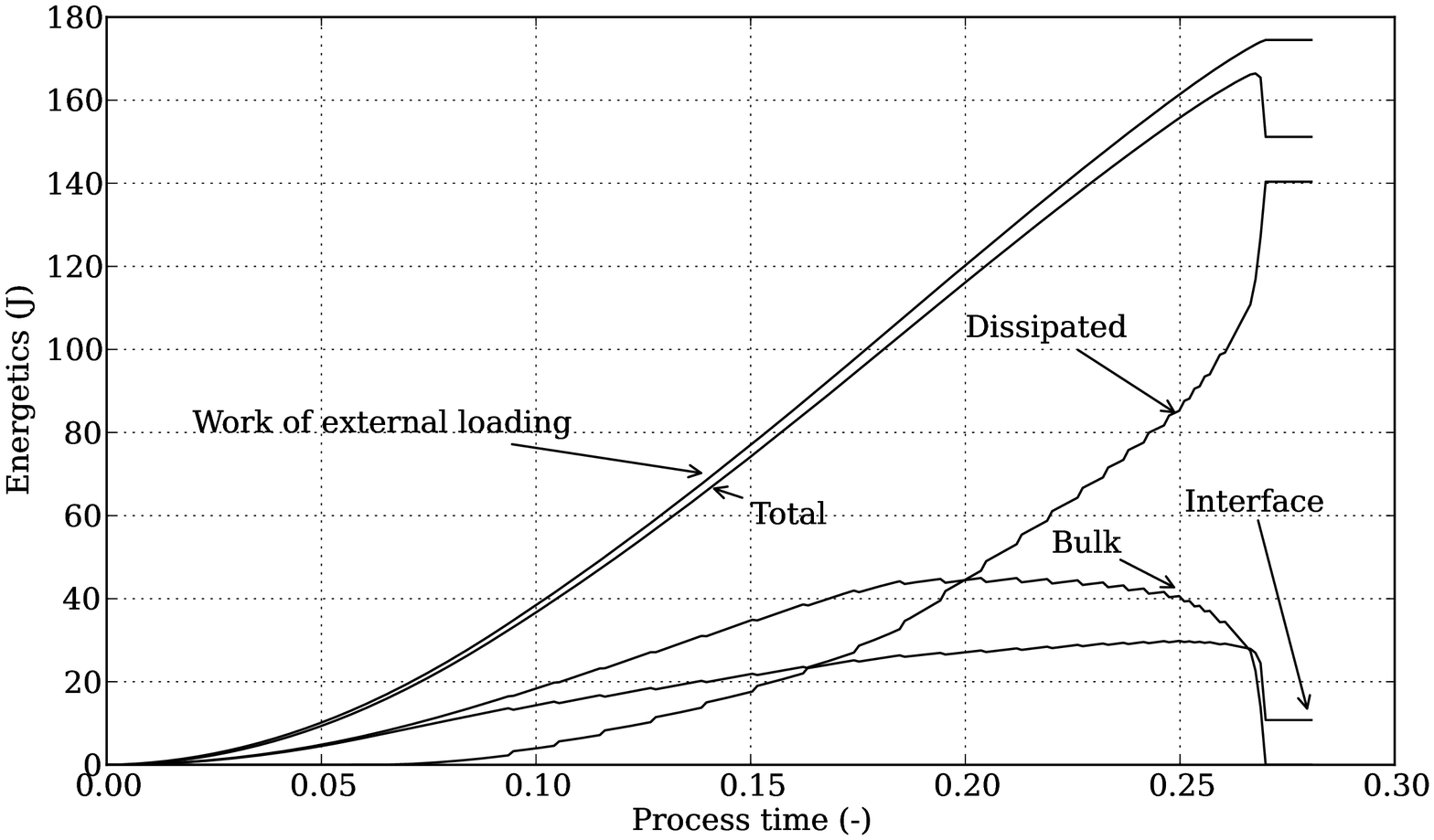}}
\end{my-picture}
\vspace*{-1em}
\begin{center}
{\small Fig.\,\ref{fig_m2}.\ }
\begin{minipage}[t]{.88\textwidth}\baselineskip=8pt
{\small
Time evolution of the energies:
the bulk and the interfacial parts of  the stored energy 
$\calE(t,\bar u_{\tau h}(t),\bar\zeta_{\tau h}(t),\bar\pi_{\tau h}(t))$,
the dissipated energy $\calR_0(\bar\zeta_{\tau h}(t){-}\zeta_0)
+\Diss_{\calR_1}^{}(\bar\pi_{\tau h};[0,t])$, their sum\,=\,total energy
(i.e.\ the left-hand side of \eqref{e6:delam-ls-mix-semi-engr}),\\[-.3em]
and the complementary work of external loading
$\int_0^t\langle\DT f_1,\underline u_{\tau h}\rangle\,\d t$
(i.e.\ the right-hand side of \eqref{e6:delam-ls-mix-semi-engr}).}
\end{minipage}
\end{center}

The bulk material is considered isotropic homogeneous with the Young modulus
$E=70$ GPa and Poisson's ratio $\nu=0.35$
(which corresponds to aluminum);
thus $\mathbb C_{ijkl}=\frac{\nu E}{(1{+}\nu)(1{-}2\nu)}\delta_{ij}\delta_{kl}
+\frac E{2(1{+}\nu)}(\delta_{ik}\delta_{jl}+\delta_{il}\delta_{jk})$
with $\delta_{ij}$ standing for the Kronecker symbol. For the adhesive, we took
a normal stiffness $\kappa_{_{\rm N}}=$150 GPa/m, a tangential
stiffness with $\kappa_{_{\rm T}}=\kappa_{_{\rm N}}/2$, the hardening slope
$\kappa_{_{\rm H}}=\kappa_{_{\rm T}}/9$, and the Mode-I fracture toughness 
$a_{_{\rm I}}=187.5$ J/m$^2$. The condition \eqref{2-sided-condition} here
means $2.65$\,MPa$\,<\sigma_{\mathrm{yield}}<5.3$\,MPa and is indeed satisfied
since $\sigma_{\mathrm{yield}}=0.56\sqrt{2 \kappa_{_{\rm N}} a_{_{\rm I}}}
=0.56\sqrt{4\kappa_{_{\rm T}}a_{_{\rm I}}}\cong4.2\,$MPa. This yields
$a_{_{\rm II}}\cong a_{_{\rm I}}+629.1\,$J/m$^2\cong816.6\,$J/m$^2$,
the fracture-mode sensitivity $a_{_{\rm II}}/a_{_{\rm I}} \cong 4.36$;
cf.\ \cite{RoMaPa13QMMD} for details.
The initial conditions are, of course, $\zeta_0=1$ and $\pi_0=0$;
the store energy $\calE\big(0,u_0,\zeta_0,\pi_0\big)$ is then 0.

It is interesting to check the energy (im)balance
\eqref{e6:delam-ls-mix-semi-engr}.
In Figure~\ref{fig_m2}, we can see it depicted for $t_1=0$ as a function of
time $t_2$: the upper line is the right-hand side of
\eqref{e6:delam-ls-mix-semi-engr} while the line below is the left-hand side
of \eqref{e6:delam-ls-mix-semi-engr}. 
We can clearly see that the difference
is not zero and is increasing in time, which is in accord with
\eqref{e6:delam-ls-mix-semi-engr} because otherwise,
if the difference would decrease on some time interval $[t_1,t_2]$,
\eqref{e6:delam-ls-mix-semi-engr} could not be valid on this interval.
This non-vanishing difference between the left- and the right-hand sides of
\eqref{e6:delam-ls-mix-semi-engr} has, beside a possible numerical error, a 
physical meaning that some part of energy is lost (dissipated) due to 
rate-dependent mechanisms, which
are neglected in the rate-independent model, like viscosity in the bulk,
cf.\ \cite{Roub13ACVE,RoPaMA13QACV}, or/and in the adhesive. One can thus
expect that, if a (vanishing) viscosity would be considered e.g.\ in the bulk,
the defect measure arising by this mechanism (like that one calculated in
\cite{RoPaMA13QACV}) would likely have the overall energy corresponding just to
this gap. Also note that, after the complete delamination, the stored energy in
the adhesive (interface) does not vanish due to the energy deposited into the
hardening.

\begin{my-picture}{.95}{.57}{fig_m3}
\begin{tabular}{lll}
\LARGE $^{^{^{^{\mbox{\footnotesize $k$=60}}}}}$  & \hspace*{-1em}\vspace*{-.1em}{\includegraphics[width=.66\textwidth]{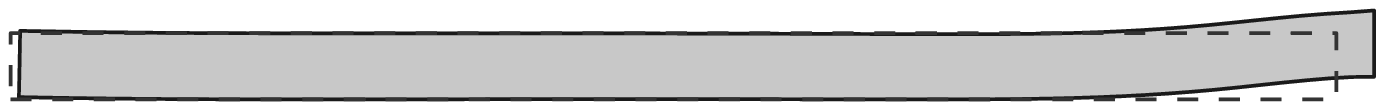}}&
\hspace*{-.5em}{\includegraphics[width=.22\textwidth,height=.13\textwidth]{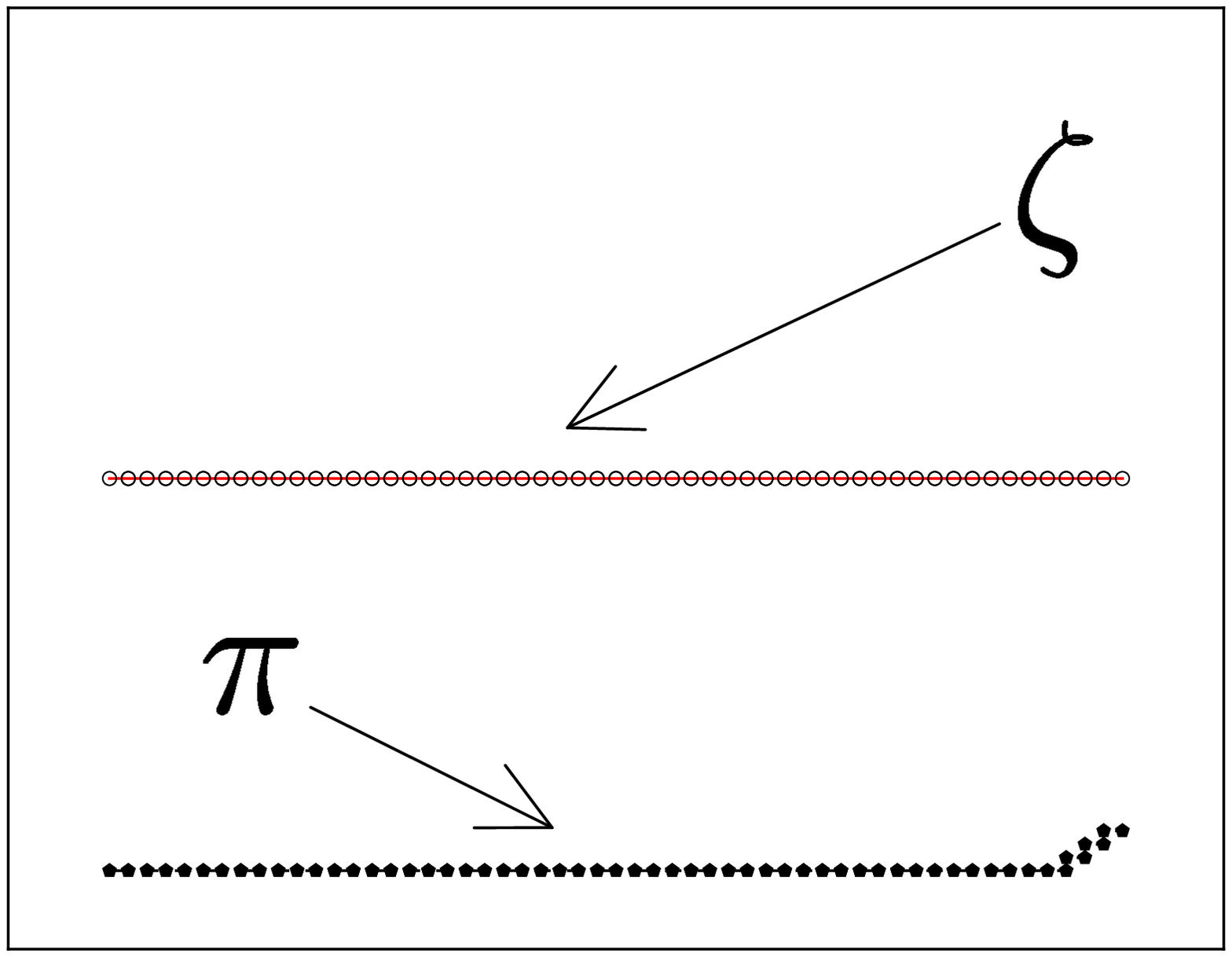}}
\\
\LARGE $^{^{^{^{\mbox{\footnotesize $k$=110}}}}}$  & \hspace*{-1em}\vspace*{-.1em}{\includegraphics[width=.66\textwidth]{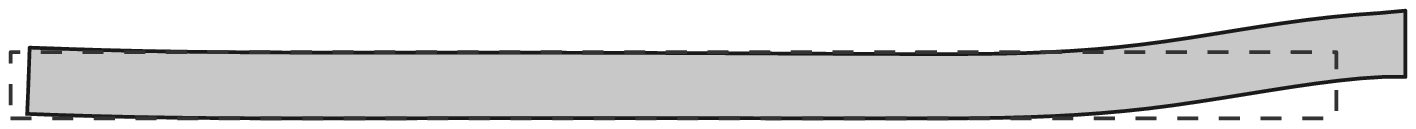}}&
\hspace*{-.5em}{\includegraphics[width=.22\textwidth,height=.13\textwidth]{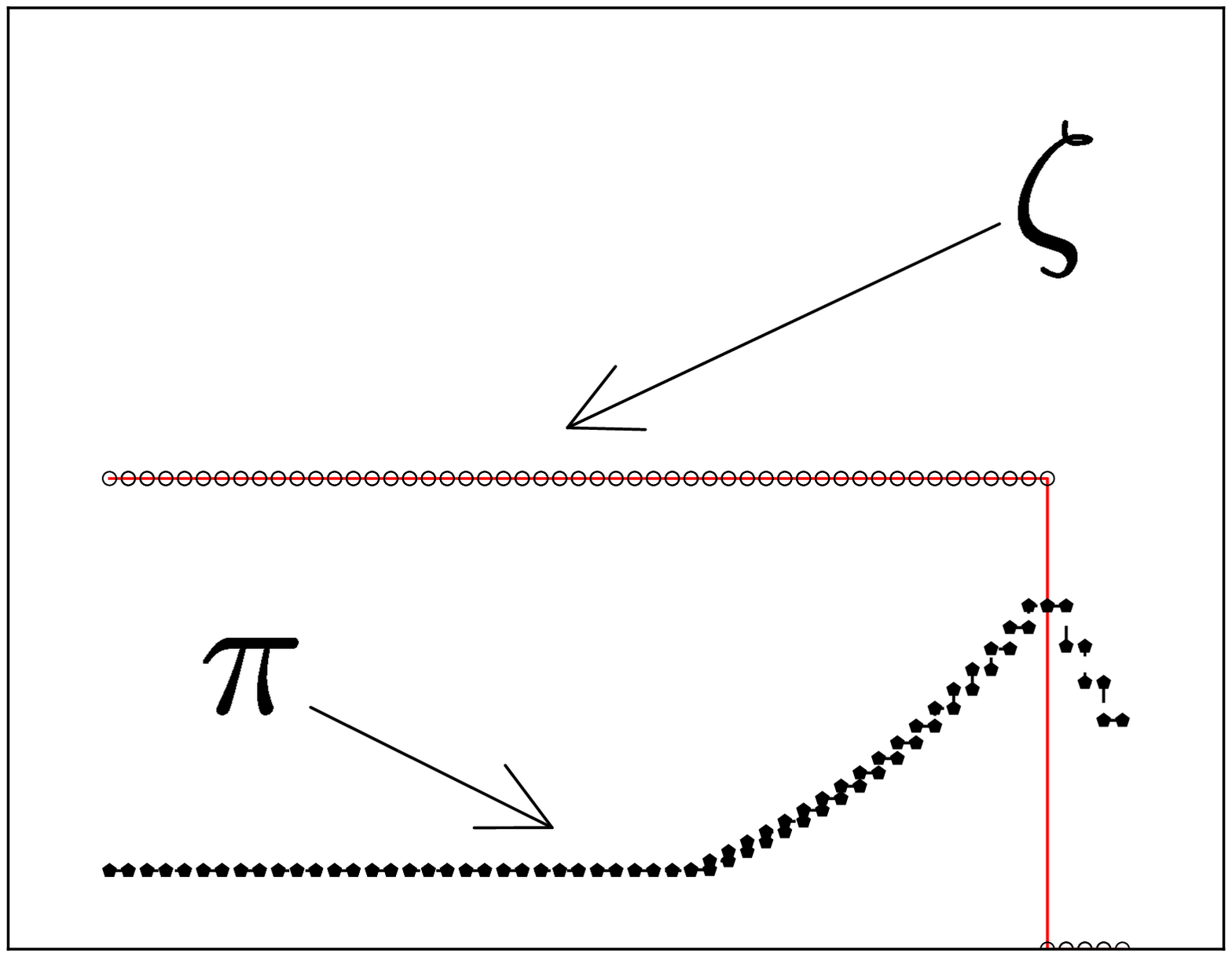}}
\\
\LARGE $^{^{^{^{\mbox{\footnotesize $k$=170}}}}}$  & \hspace*{-1em}\vspace*{-.1em}{\includegraphics[width=.66\textwidth]{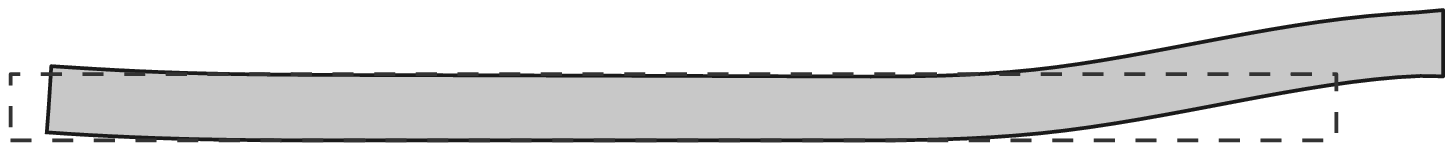}}&
\hspace*{-.5em}{\includegraphics[width=.22\textwidth,height=.13\textwidth]{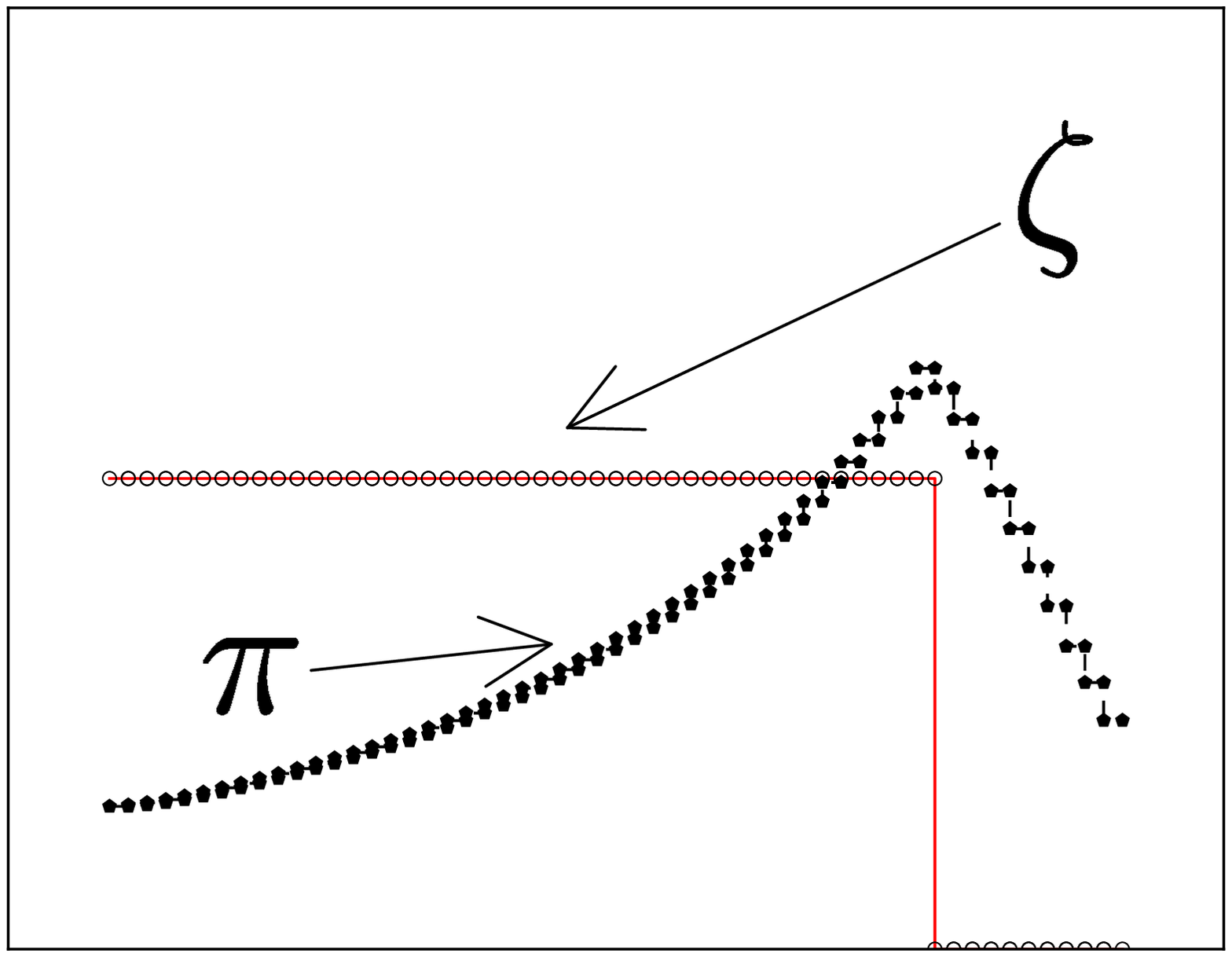}}
\\
\LARGE $^{^{^{^{\mbox{\footnotesize $k$=210}}}}}$  & \hspace*{-1em}\vspace*{-.1em}{\includegraphics[width=.66\textwidth]{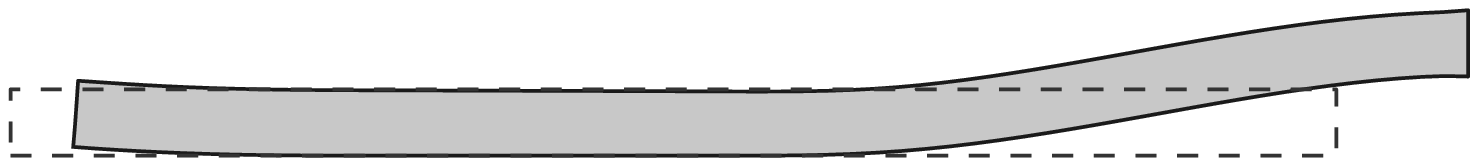}}&
\hspace*{-.5em}{\includegraphics[width=.22\textwidth,height=.13\textwidth]{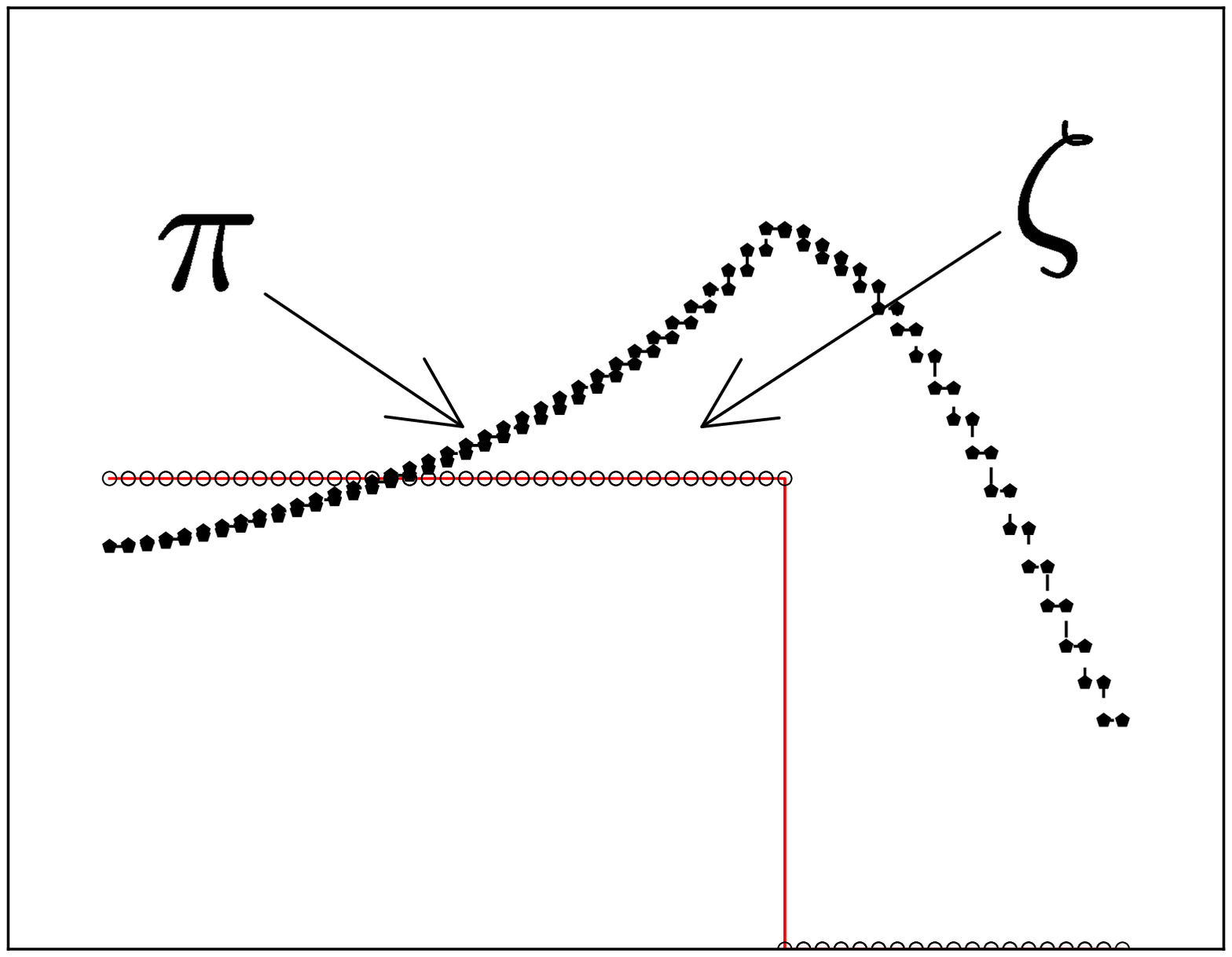}}
\\
\LARGE $^{^{^{^{\mbox{\footnotesize $k$=225}}}}}$  & \hspace*{-1em}\vspace*{-.1em}{\includegraphics[width=.66\textwidth]{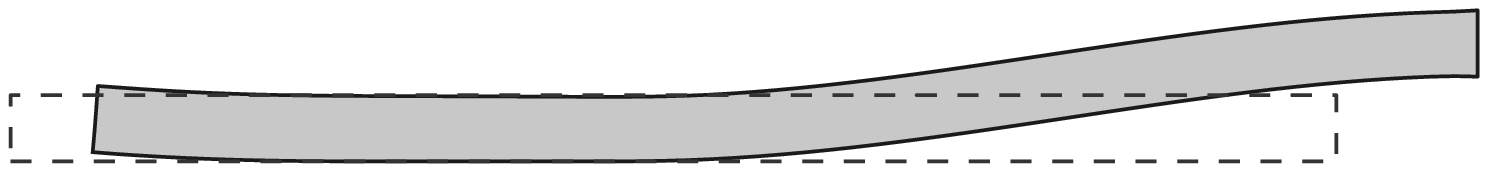}}&
\hspace*{-.5em}{\includegraphics[width=.22\textwidth,height=.13\textwidth]{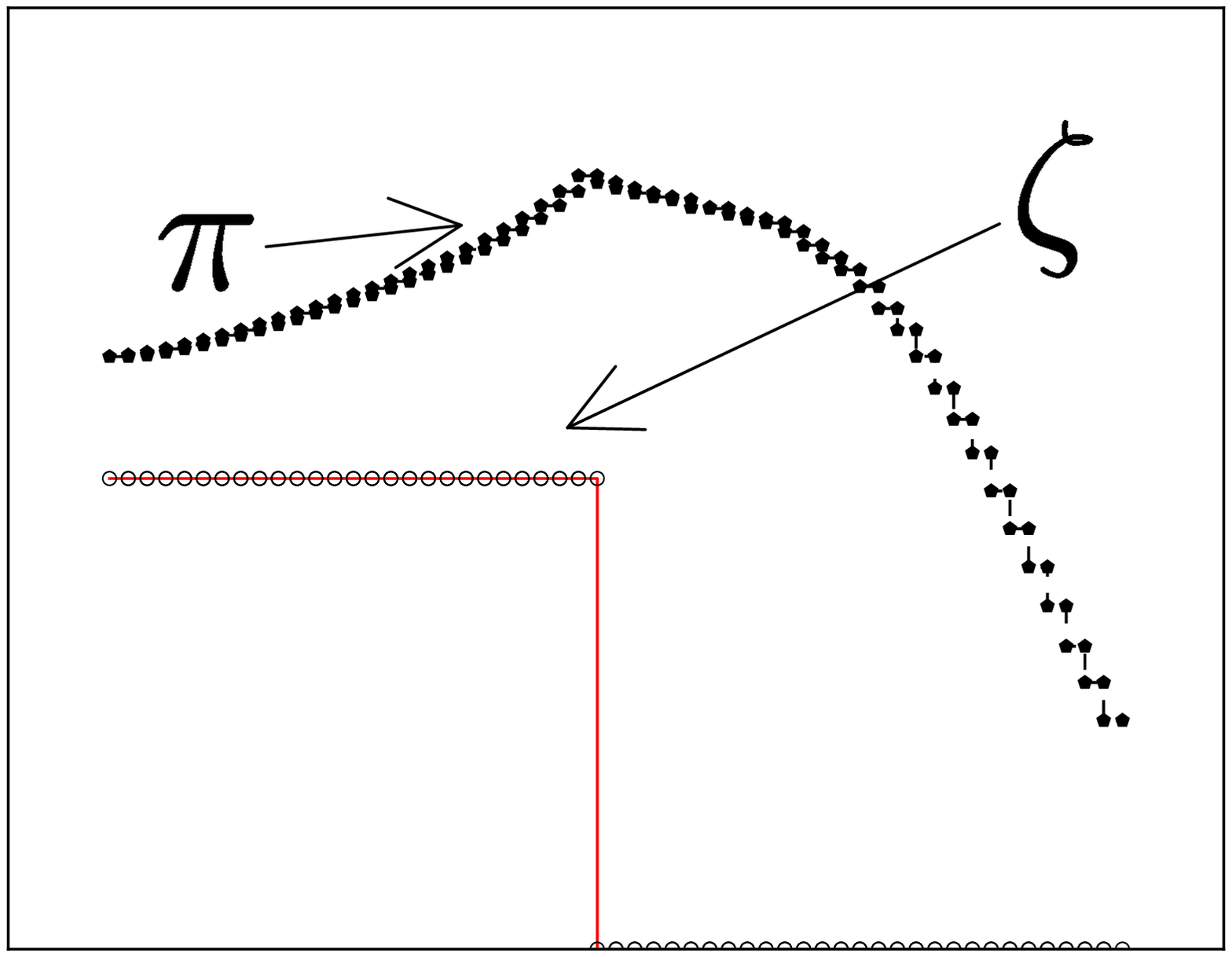}}
\\
\LARGE $^{^{^{^{\mbox{\footnotesize $k$=226}}}}}$  & \hspace*{-1em}\vspace*{-.1em}{\includegraphics[width=.66\textwidth]{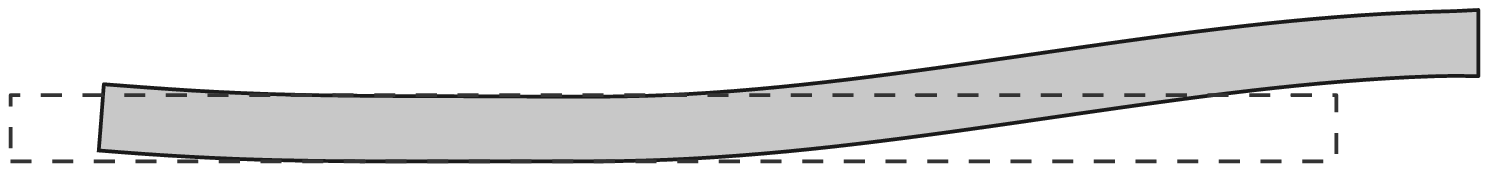}}&
\hspace*{-.5em}{\includegraphics[width=.22\textwidth,height=.13\textwidth]{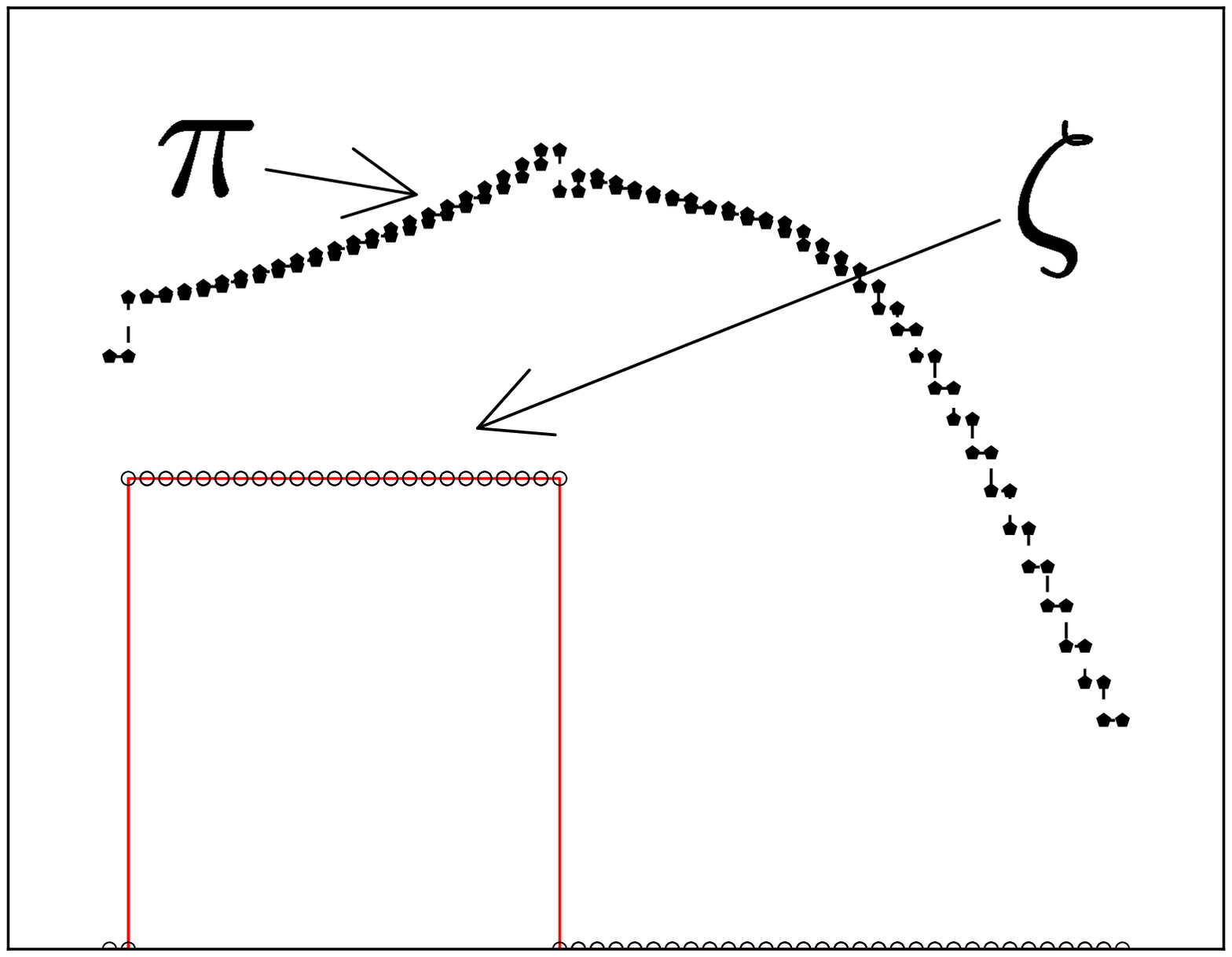}}
\\
\LARGE $^{^{^{^{\mbox{\footnotesize $k$=227}}}}}$  & \hspace*{-1em}\vspace*{-.1em}{\includegraphics[width=.66\textwidth]{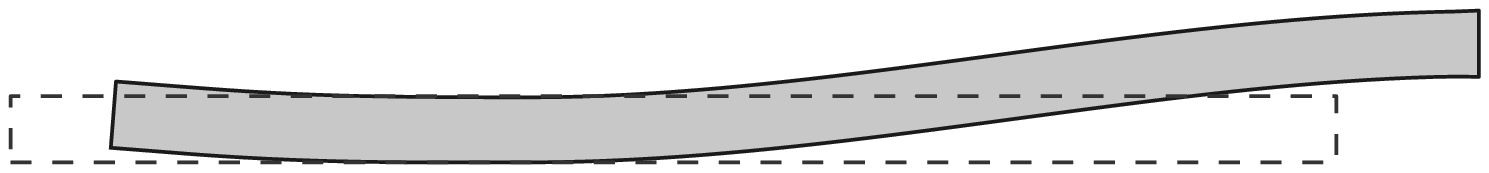}}&
\hspace*{-.5em}{\includegraphics[width=.22\textwidth,height=.13\textwidth]{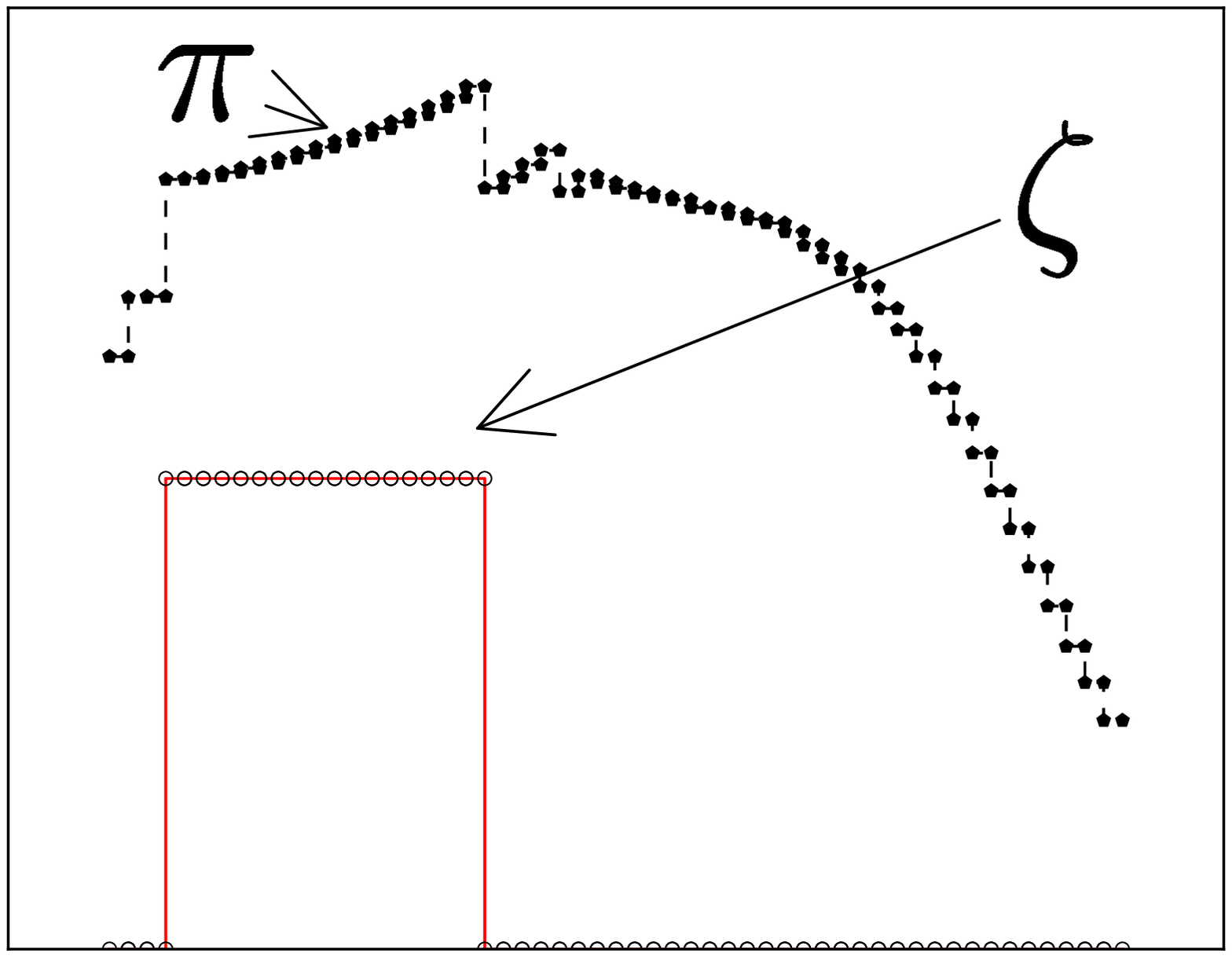}}
\\
\LARGE $^{^{^{^{\mbox{\footnotesize $k$=228}}}}}$  & \hspace*{-1em}\vspace*{-.1em}{\includegraphics[width=.66\textwidth]{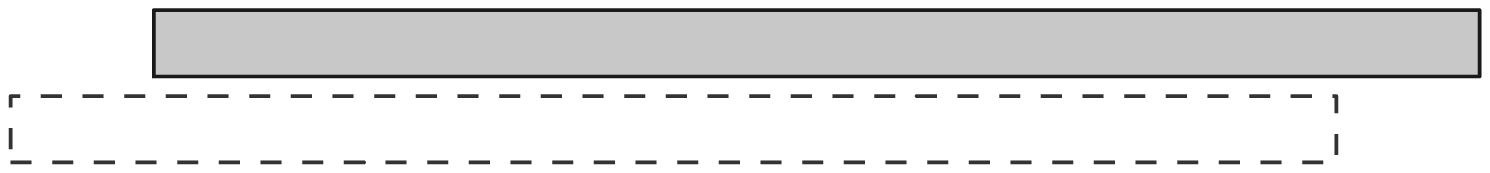}}&
\hspace*{-.5em}{\includegraphics[width=.22\textwidth,height=.13\textwidth]{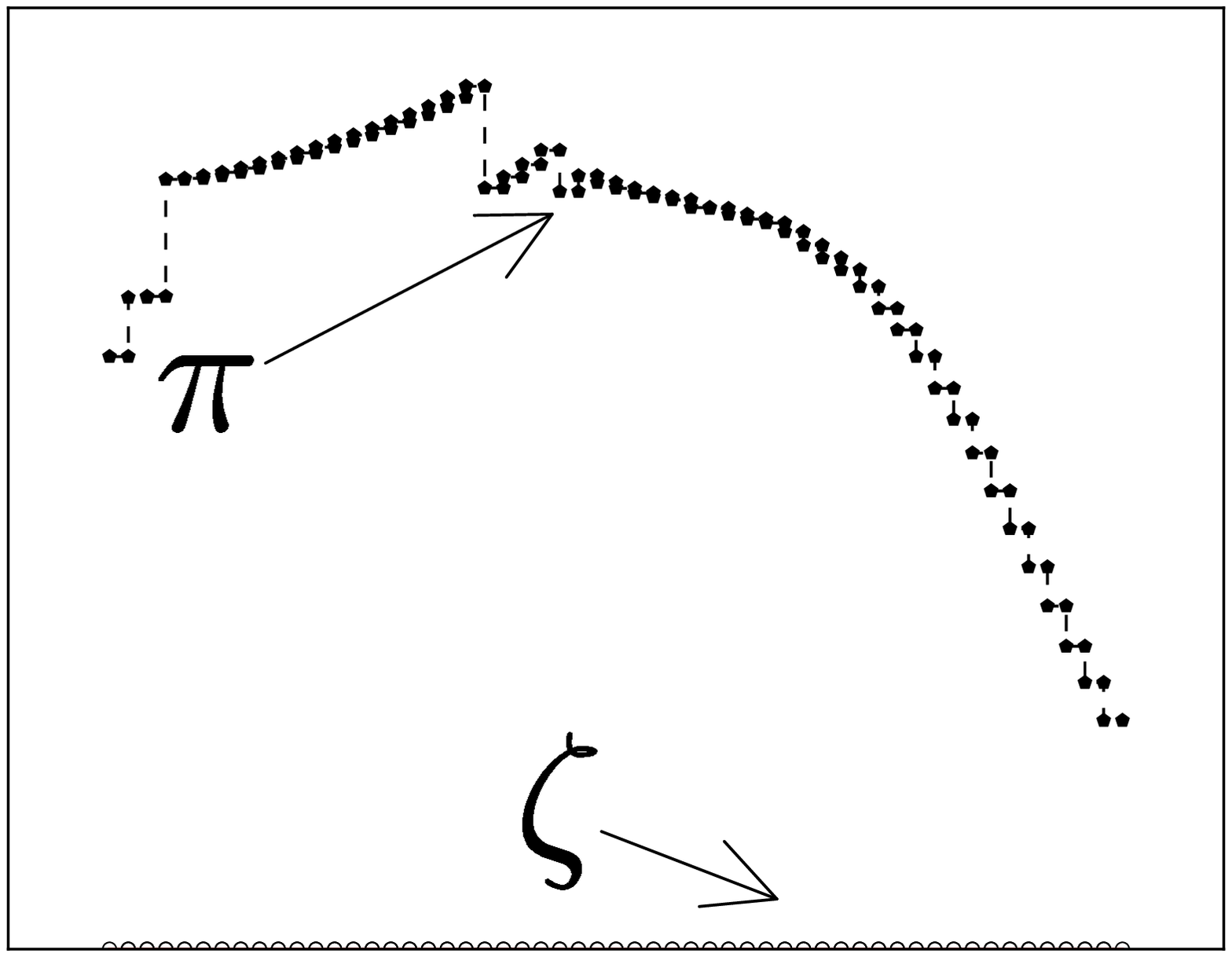}}
\end{tabular}
\end{my-picture}
\vspace*{19em}
\begin{center}
{\small Fig.\,\ref{fig_m3}.\,}
\begin{minipage}[t]{.89\textwidth}\baselineskip=8pt
{\small
Time evolution at eight snapshots of the geometrical configuration
(displacement depicted magnified $100\,\times$) and the spatial distribution of
$\zeta$ and $\pi$ along $\GC$.
}
\end{minipage}
\end{center}

This example exhibits remarkably varying mode of delamination. At the beginning
the delamination is performed by a mixed mode close to Mode I given
essentially by the direction of the Dirichlet loading, cf.\ Figure~\ref{fig_m1},
while later it turns rather to nearly pure Mode II. Yet, at the very end of 
the process, due to elastic bending the delamination starts performing also 
from the left-hand side of the bar opposite to the loading side, and thus 
again a mixed mode occurs. This relatively complicated mixed-mode behaviour
is depicted in Figures~\ref{fig_m3}--\ref{fig_m4}, showing essential 
qualitative difference from the energetic solution which exhibits a 
non-physical tendency to slide to less-dissipative Mode I, 
cf.\ \cite[Fig.\,7]{RoMaPa13QMMD}.
%\vspace*{-17em}\vspace*{-0em}

The evolution of the deformation $u$ and spatial distribution of
the delamination $\zeta$ and the plastic slip $\pi$ are depicted in
Figure~\ref{fig_m3} at eight snapshots selected not uniformly to
visualize interesting effects when delamination starts to be completed.
In particular, the delamination propagating from both
sides at the very end (mentioned already above) is seen there.

\begin{my-picture}{.95}{.35}{fig_m4}
\hspace*{.5em}{\includegraphics[width=.47\textwidth]{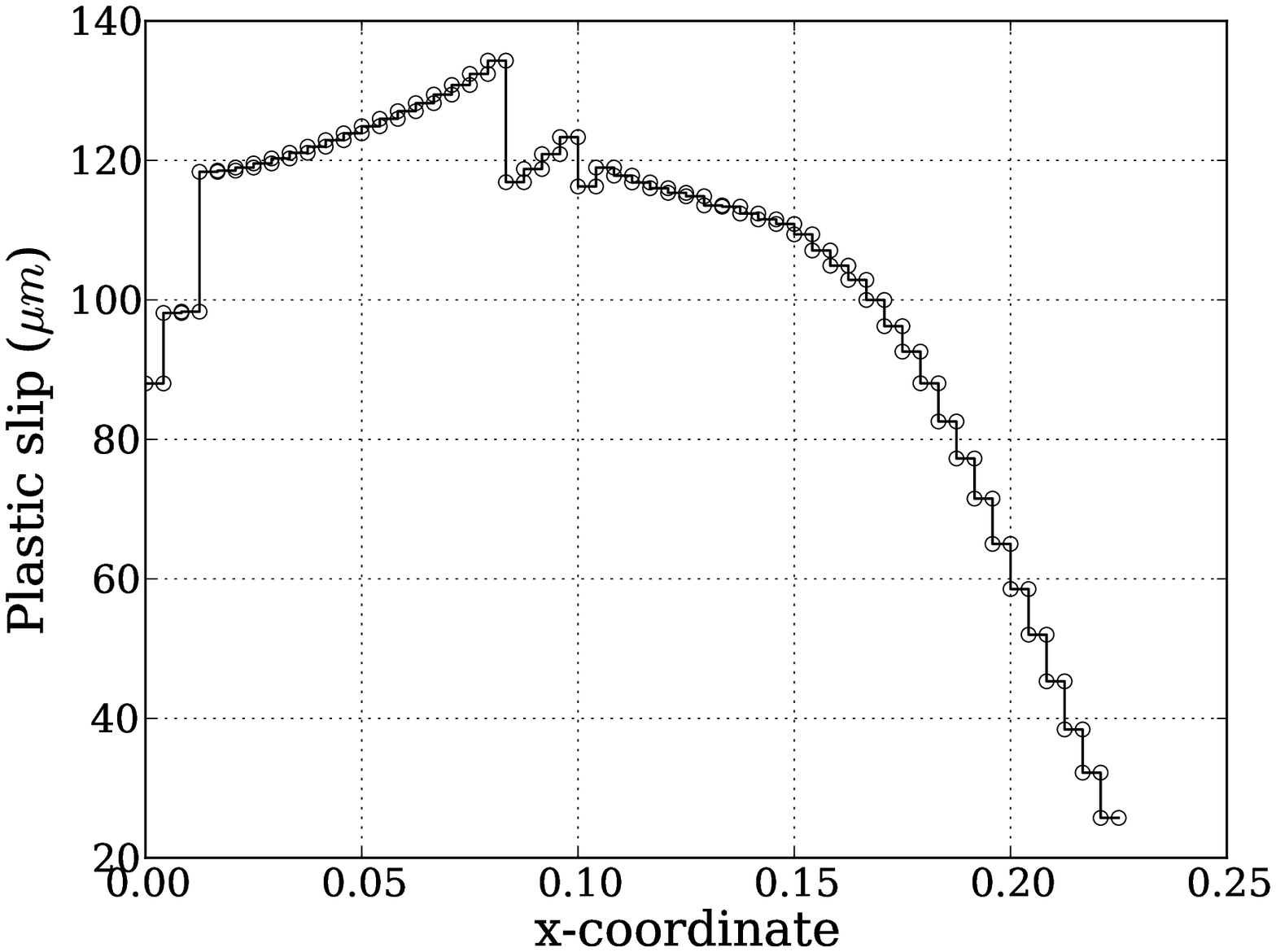}}
\hspace*{.5em}\vspace*{-.1em}{\includegraphics[width=.47\textwidth]{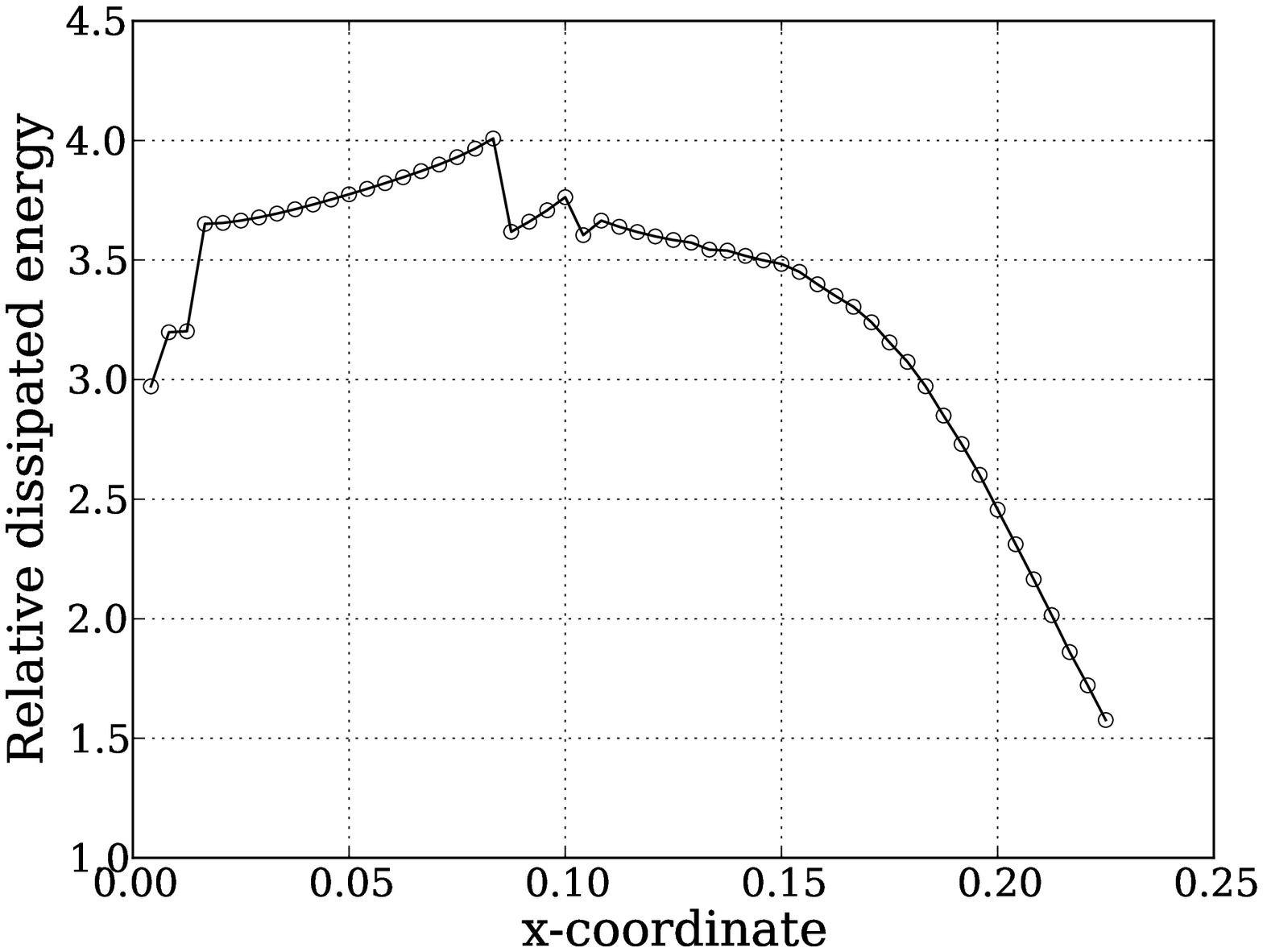}}
\end{my-picture}
\vspace*{-1em}
\begin{center}
{\small Fig.\,\ref{fig_m4}.}
\begin{minipage}[t]{.91\textwidth}\baselineskip=8pt
{\small 
Distribution of mode-mixity of delamination  along $\GC$:
\\Left:
The overall plastic slip after the delamination has been completed
(=the last snapshot in Figure~\ref{fig_m3}).
\\
Right: The dissipated energy related to $a_{_{\rm I}}$ after the delamination 
has been completed (value=1$\sim$Mode I, value=3.97 $\sim$Mode II).
Similar distributions are observed in both plots  because there were not 
cycling in plastification during the delamination.
}
\end{minipage}
\end{center}

For the discretisation of the experiment in Figures~\ref{fig_m2}--\ref{fig_m4},
we choose $\tau=0.012$ and $h=4.6\,$mm (=the size of a boundary element in
uniform  discretisation).

\begin{my-picture}{.95}{.4}{fig-AMDP}
\hspace*{-1em}\vspace*{-0em}{\includegraphics[width=.52\textwidth]{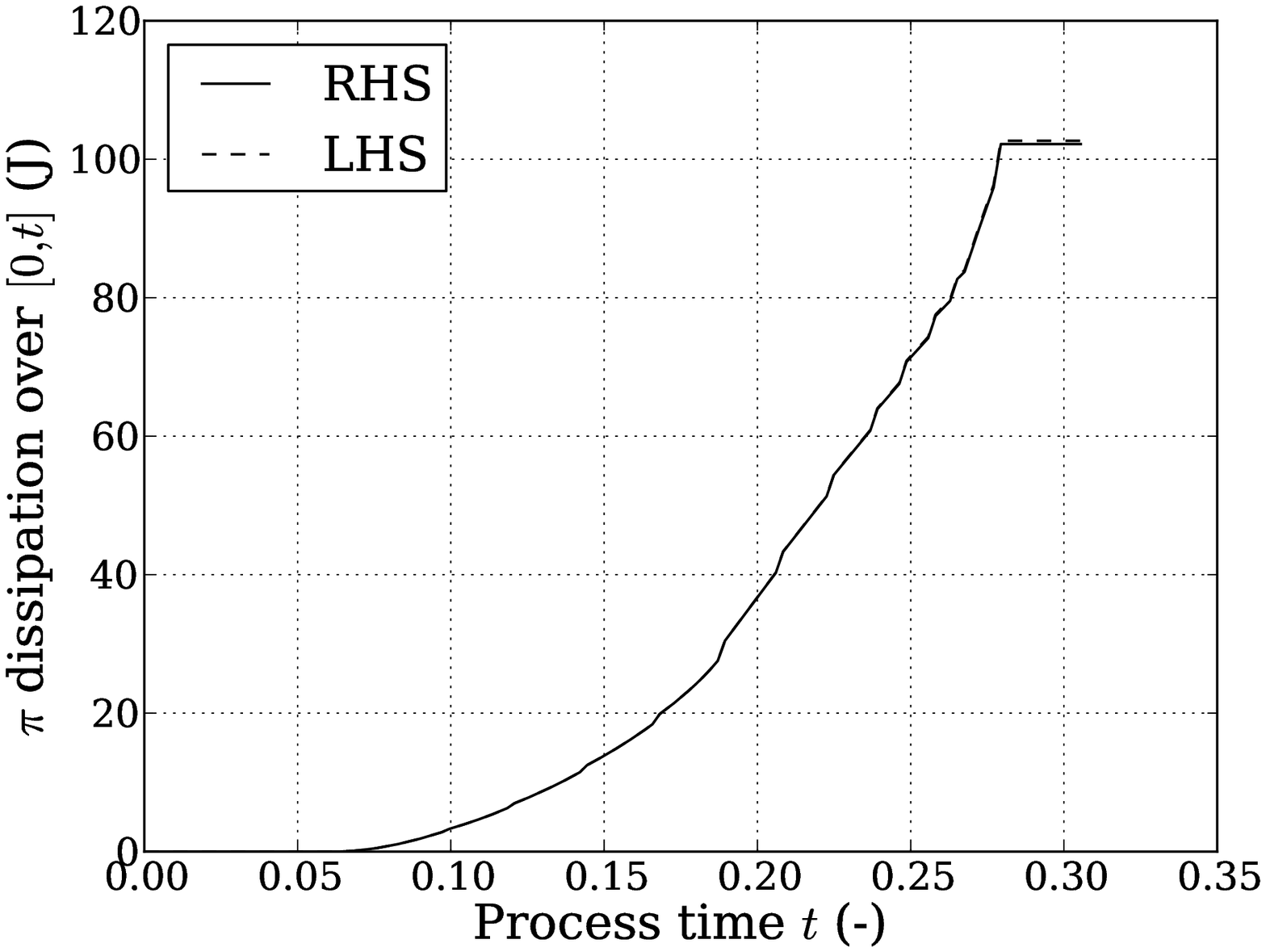}}
\hspace*{-1em}\vspace*{-0em}{\includegraphics[width=.52\textwidth]{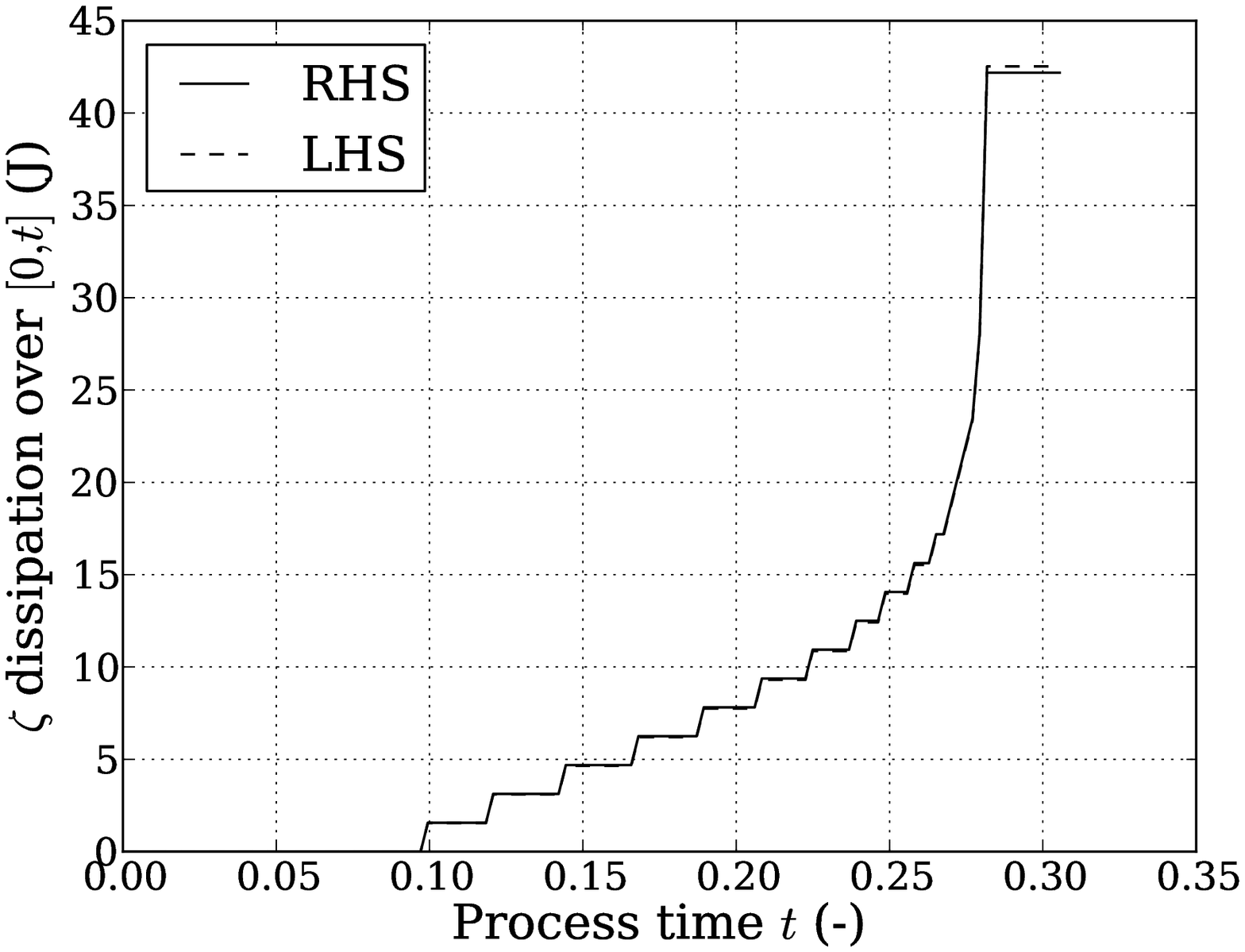}}
\end{my-picture}
\vspace*{-1em}
\begin{center}
{\small Fig.\,\ref{fig-AMDP}.\ }
\begin{minipage}[t]{.8\textwidth}\baselineskip=8pt
{\small
Time evolution of the left- and the right-hand sides in
the approximate maximum dissipation principle (AMDP)
for the plastic slip $\pi$, i.e.\ \eqref{e1:AMDP+a},
and the damage parameter $\zeta$, i.e.\ \eqref{e1:AMDP+b}. The difference
is practically invisible in the former case
and less than 2\% in the latter case.
}
\end{minipage}
\end{center}
The differences in the approximate maximum-dissipation principles 
\eqref{e1:AMDP+} are now displayed in Figure~\ref{fig-AMDP}. We can see that 
our algorithm yielded a well (about 98\%) maximally-dissipative (i.e.\ 
stress-driven) solution, the possible deviation is possibly only in $\zeta$ 
at the very end of the delamination process. 

Eventually, the joint convergence from Corollary~\ref{cor-convergence}
for time- and FEM-spatial discretisation (although here implemented by BEM) 
is demonstrated in Figures~\ref{fig_conv} and \ref{fig_conv+} for a twice 
coarser time/space discretisations. We choose the scenario keeping the ratio 
$\tau/h$ constant, although Corollary~\ref{cor-convergence} itself does not 
give any particular suggestion in this respect. Anyhow, the tendency of
convergences is clearly seen, although we naturally do not know the
exact solution so that we cannot evaluate any actual error. On top of it,
the exact solution does not need to be unique so we even do not have guaranteed
the convergence of the whole sequence of the approximate solutions and, 
moreover, the simplified implementation by
collocation BEM does not have guaranteed convergence,
in contrast to FEM stated in Corollary~\ref{cor-convergence}.

\begin{my-picture}{.95}{.35}{fig_conv}
\hspace*{-.5em}\vspace*{-.1em}{\includegraphics[width=.45\textwidth]{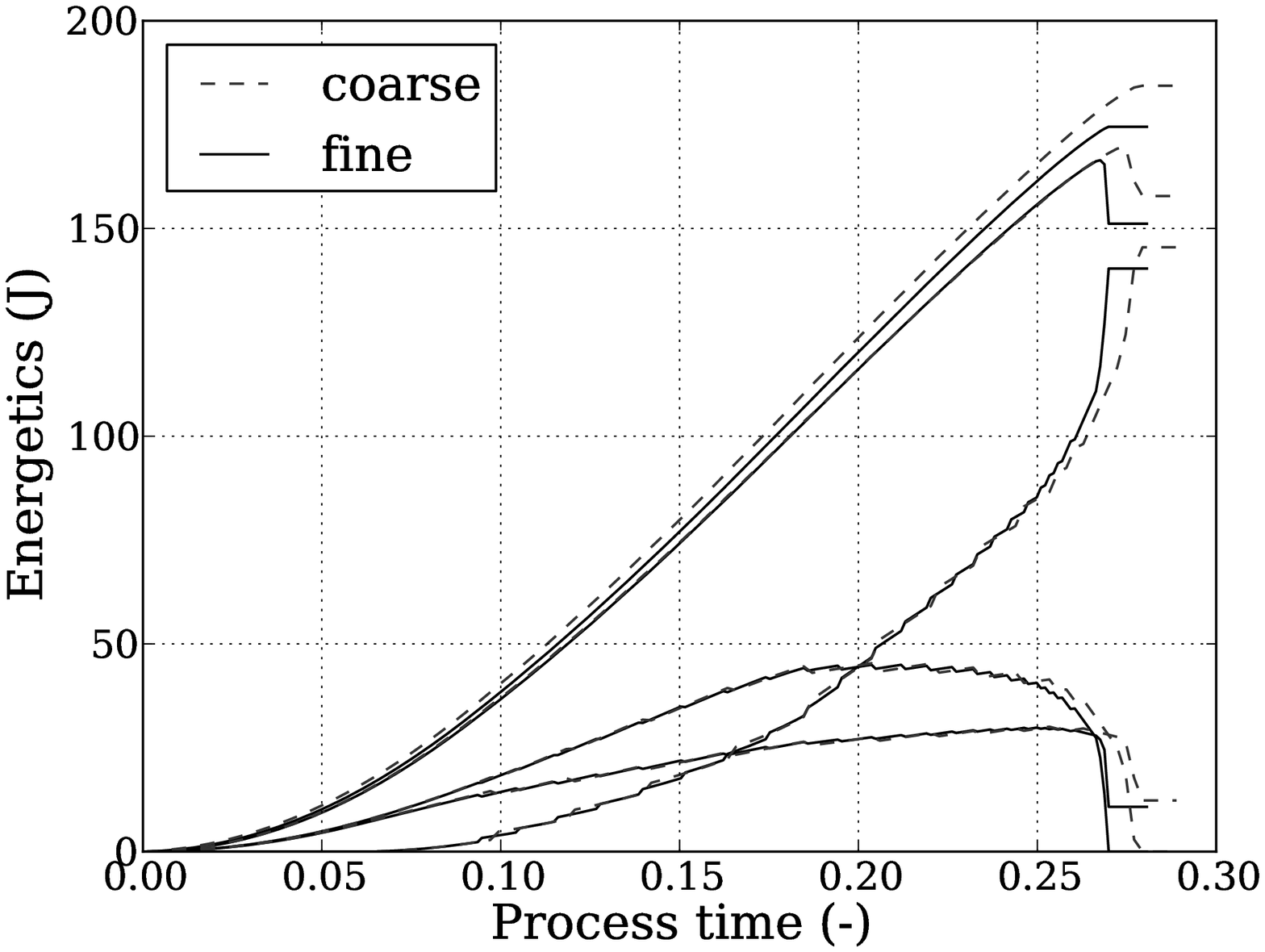}}
\hspace*{.5em}\vspace*{-.1em}{\includegraphics[width=.45\textwidth]{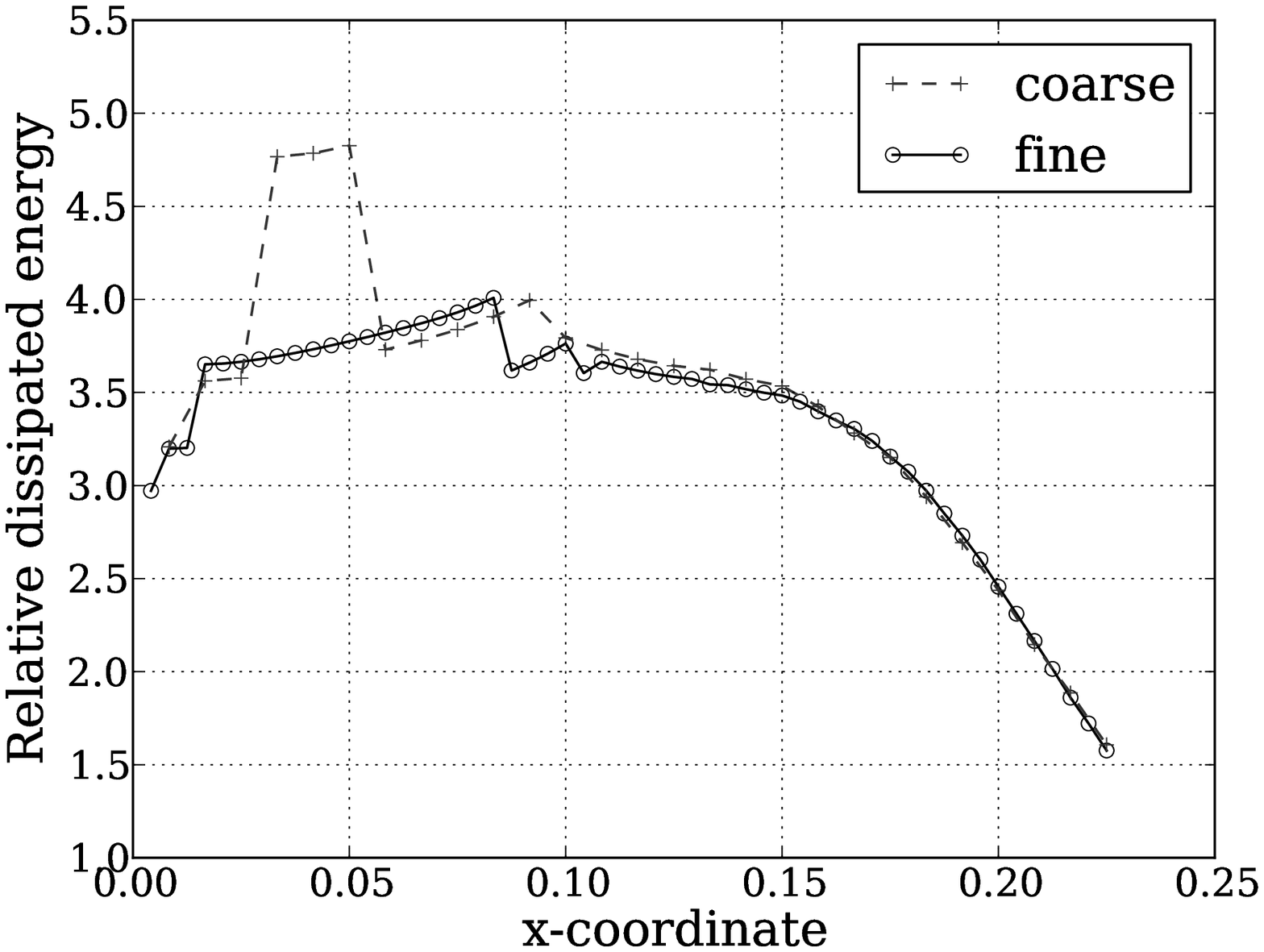}}
\end{my-picture}
\vspace{-1em}
\begin{center}
{\small Fig.\,\ref{fig_conv}.\ }
\begin{minipage}[t]{.8\textwidth}\baselineskip=8pt
{\small
Convergence test: Left:\ \ evolution of energies as in Figure~\ref{fig_m2}
\\Right: final spatial distribution  of $\pi$ along $\GC$ as in Figure~\ref{fig_m4}(right).
}
\end{minipage}
\end{center}

\begin{my-picture}{.95}{.35}{fig_conv+}
\hspace*{-.5em}\vspace*{-.1em}{\includegraphics[width=.45\textwidth]{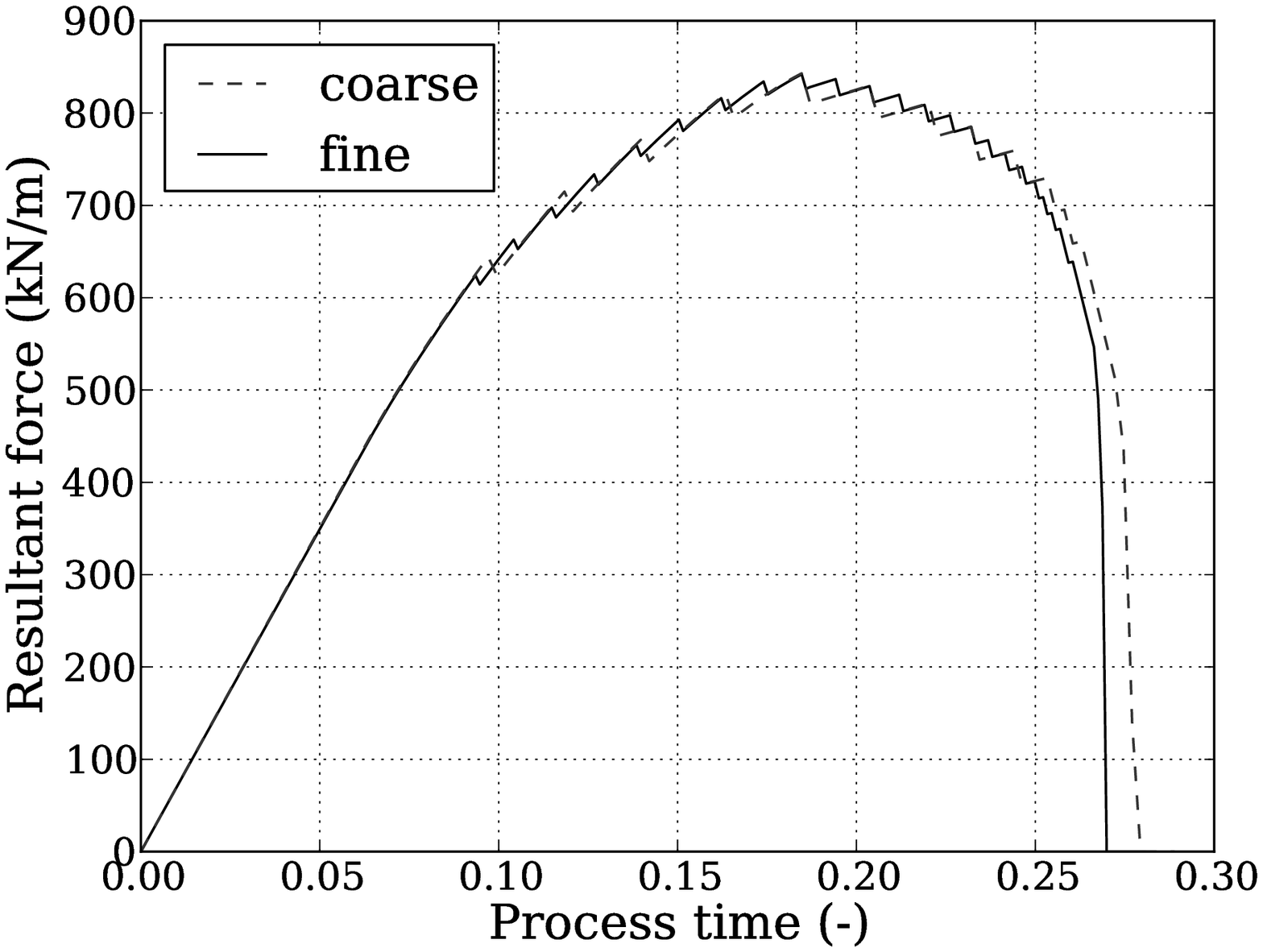}}
\hspace*{.5em}\vspace*{-.1em}{\includegraphics[width=.45\textwidth]{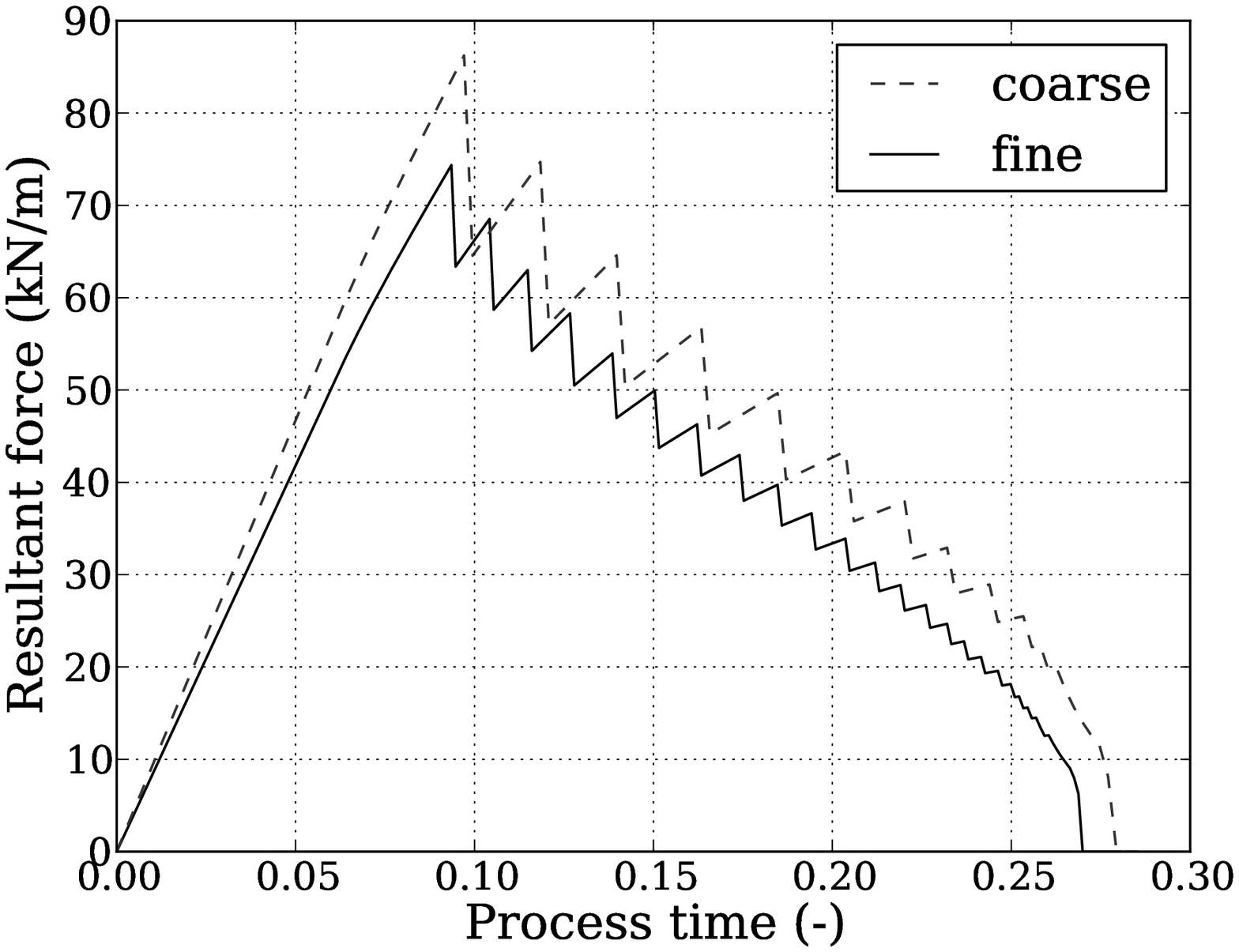}}
\end{my-picture}
\vspace{-1em}
\begin{center}
{\small Fig.\,\ref{fig_conv+}.\ }
\begin{minipage}[t]{.8\textwidth}\baselineskip=8pt
{\small
Convergence test: horizontal (left) and vertical (right) component of the 
total force response evolving in time.
}
\end{minipage}
\end{center}

\baselineskip=11pt

\subsection*{Acknowledgments}
The authors are thankful to anonymous referees for many valuable comments
that led to improvement of presentation in many spots. This research has been 
covered by the Junta de Andaluc\'{\i}a (Proyecto de Excelencia P08-TEP-4051) 
and  the Spanish Ministry of Economy and Competitiveness (MAT2012-37387) as 
well as partial support from the grants 201/10/0357 and 13-18652S
(GA \v CR), together with the institutional support RVO:\,61388998 (\v CR).
T.R.\  (resp.\ C.G.P.)  acknowledges the hospitality of
Universidad de Sevilla, where this work has partly (resp.\ mostly)
been accomplished.

\bibliographystyle{abbrv}

\bibliography{trsevilla7}

\end{document}